\documentclass[sts]{imsart}

\RequirePackage[colorlinks,citecolor=blue,urlcolor=blue]{hyperref}

\RequirePackage{amsthm,amsmath,amsfonts,amssymb}
\RequirePackage[authoryear]{natbib}
\RequirePackage[colorlinks,citecolor=blue,urlcolor=blue]{hyperref}
\RequirePackage{graphicx}

\usepackage{enumitem}
\usepackage{url}
\usepackage[super]{nth}
\usepackage[utf8]{inputenc}
\usepackage[normalem]{ulem}
\usepackage{scalerel}
\usepackage{tabularx}
\usepackage{subfig}
\usepackage{dsfont}

\newcommand{\PageRank}{\mbox{PageRank}}

\begin{document}

\begin{frontmatter}
\title{The many routes to the ubiquitous Bradley--Terry model}
%\title{A Sample Article Title with Some Additional Note\thanksref{t1}}
\runtitle{The many routes to the ubiquitous Bradley--Terry model}
%\thankstext{T1}{A sample additional note to the title.}

\begin{aug}
\author[A]{\fnms{Ian}~\snm{Hamilton}\ead[label=e1]{ian.hamilton@warwick.ac.uk}\orcid{0000-0002-3991-9693}},
\author[B]{\fnms{Nicholas}~\snm{Tawn}\ead[label=e2]{n.tawn.1@warwick.ac.uk}}
\and
\author[C]{\fnms{David}~\snm{Firth}\ead[label=e3]{d.firth@warwick.ac.uk}\orcid{0000-0003-0302-2312}}

\address[A]{Ian Hamilton is Honorary Fellow, Department of Statistics,
University of Warwick, Coventry, U.K.\printead[presep={\ }]{e1}.}

\address[B]{Nicholas Tawn is Associate Professor of Statistics, Department of Statistics,
University of Warwick, Coventry, U.K.\printead[presep={\ }]{e2}.}

\address[C]{David Firth is Emeritus Professor of Statistics, Department of Statistics,
University of Warwick, Coventry, U.K.\printead[presep={\ }]{e3}.}

\end{aug}

\begin{abstract}
The rating of items based on pairwise comparisons has been a topic of statistical investigation for many decades. Numerous approaches have been proposed. One of the best known is the Bradley--Terry model. This paper seeks to assemble and explain a variety of motivations for its use. Some are based on principles or on maximizing an objective function; others are derived from well-known statistical models, or stylized game scenarios. They include both examples well-known in the literature as well as what are believed to be novel presentations.
\end{abstract}


\end{frontmatter}

%\addtocontents{toc}{\protect\setcounter{tocdepth}{-1}}

\section{Introduction}
The first conference that the lead author attended as a PhD student was an American sports statistics conference. He presented a poster related to the Bradley--Terry model. As a retrodictive model on rugby union in a sea of American sports predictions, it felt a little out of place. But a kind attendee took pity on him and decided to engage him with a question. She asked, ``Why would I choose Bradley--Terry rather than the Thurstone model?" (by which he took her to mean what is more commonly referred to as the Thurstone--Mosteller model). He flummered a vague response involving analytic niceness and simplicity --- he suspects Occam's razor even got a mention. She looked suitably unconvinced. It is to be hoped that this paper represents a more ordered response to the conference interlocutor and an aggregation of, as \citet[p.13]{david1963method} puts it in his canonical survey of pairwise comparison methods, ``the many routes to the ubiquitous Bradley--Terry model.'' 

%This may be seen as an updating and expansion of Section 3 of \cite{bradley1976science}, with the presentations of Sections \ref{sec: objective function maximisation}, \ref{sec: standard models}, \ref{sec: games}, and \ref{sec: quasi-symmetry} here not discussed in that work and those of Sections \ref{sec: axiomatic} and \ref{sec: discriminal processes} expanded upon.

Thus, the main original contribution of the work is in aggregating the motivations for the Bradley--Terry model, or as \citet{bradley1976science} refers to them, the `bases for model formulation'. In collating these motivations, we hope that the work provides a useful resource to those encountering the model for the first time and some new perspectives for those more familiar with it. It may also complement other works, such as \citet{david1963method}, \citet{cattelan2012models}, \citet{vojnovic2015contest}, or \citet{wu2022diagnostic} that provide alternative helpful summary perspectives on the model. The work takes in a diverse scope of motivating ideas including likelihood and entropy maximization, psychological choice and sensation models, distance minimization, a prominent Markov chain Monte Carlo method, other well-known rating models such as \PageRank\ and the RPI of American college sports, sudden-death play-offs, pub pool norms and the British playground game of conkers. The aggregation of these motivations serves to demonstrate the broad appeal of the Bradley--Terry model in many settings. 

The paper also offers a number of novelties including: a more extensive explicit discussion of the Bradley--Terry model in the context of an exponential family of distributions than has appeared previously, which provides a uniting theme to a number of the more notable motivations; a formalization of perhaps the most intuitive motivation for the model, by proposing an explicit measure for the simplicity of a model in the pairwise comparison scenario and showing that, under plausible constraints, Bradley--Terry is the model that maximizes this measure; and a demonstration of how the ideas behind the rating method of \citet{wei1952algebraic} and \citet{kendall1955further} and of the Ratings Percentage Index (RPI) can be related to the Bradley--Terry model through the Perron--Frobenius Theorem.

The scenario under consideration in this paper is one where there is a desire to create a ranking of items based on the observation of a set of binary-outcome pairwise comparisons. One popular approach to ranking is to determine a uni-dimensional rating, and then order items by their ratings. Statistical models such as Bradley--Terry or Thurstone--Mosteller achieve this by defining the probability of a preference for alternative $i$ over alternative $j$ in a pairwise comparison independently from other preferences conditional on the strengths of the items. In the Bradley--Terry model the probability is defined as
\[
p_{ij} = \frac{\pi_i}{\pi_i + \pi_j},
\]
where $\pi_i$ is a positive-valued parameter that may be interpreted as a rating of alternative $i$, with a higher rating indicating greater `strength' or `worth'. 

This results in a model that generates independent Binomial realizations between pairs of items. Therefore, with a logit transformation of the above, one can equivalently state the model as a member of the class of generalized linear models \citep{mccullagh1989generalized} with
\[
F(p_{ij})=\lambda_i-\lambda_j ,
\]
where $\lambda_i = \log(\pi_i)$ is a real-valued parameter indicating the strength of $i$, and $F$ is taken as the logit function. The Thurstone--Mosteller model \citep{thurstone1927law, mosteller1951remarks}, about which the interlocutor asked, is derived from taking $F$ to be the probit function instead.  In practice, the fitted models are usually very similar \citep{chambers1967discrimination, stern1992all}.

%\footnote{Interesting motivations for the use of a generalised linear model itself in the context of pairwise comparison include those due to \citet{brunk1960mathematical}, in showing how it may be derived from a utility model, and \citet{joe1988majorization}, in showing that such a model is the maximum entropy model with respect to a majorisation ordering of the possible probability matrices.}

The Bradley--Terry model has formed the basis for many models and analyses in many contexts over time.  These include, for example, analysis of journal citations  \citep{stigler1994citation}, college sports \citep{wobus2007KRACH}, animal behavior \citep{stuart2006multiple}, risk analysis \citep{merrick2002prince}, wine tasting \citep{oberfeld2009ambient}, university ranking \citep{dittrich1998modelling}, font selection \citep{o2014exploratory}, educational assessment \citep{pollitt2012method}, development of large language models \citep[especially following][]{rafailov2023dpo} --- and of course chess, which was the subject of the original work by \citet{zermelo1929berechnung} as well as being the subject of the popular closely-related ranking method proposed by \citet{elo1978rating}, still widely used in chess today. Extended versions of the model include those addressing dynamic ranking \citep[e.g.,][]{Glickman1999dynamic, Cattelan2013dynamic}, item-specific and/or judge-specific covariates \citep[e.g.,][]{schauberger2019covariates}, random effects \citep[e.g.,][]{Lancaster1983randomeffects, Matthews1995randomeffects, Bockenholt2001randomeffects} and spatial proximity \citep{seymour2022bayesian}.

Originally documented by \citet{zermelo1929berechnung}, the Bradley--Terry model took the name by which it came to be commonly known when \citet{bradley1952rank} independently rediscovered it. Following the work of \citet{thurstone1927law, thurstone1927method,thurstone1927psychophysical} and \citet{zermelo1929berechnung}, paired comparison methods saw little development for the best part of a quarter of a century until they became an active area of investigation in the 1950s and 60s. Much of this work took place in the context of the psychological literature, with Luce's Choice Axiom \citep{luce1959individual} a particularly notable contribution, leading to the model sometimes being referred to as the Bradley--Terry-Luce (BTL) model. A number of these works showed how the Bradley--Terry model could be derived based on plausible axioms or desirable model features \citep{good1955marking, luce1959individual, buhlmann1963pairwise, luce1965preference}. Towards the end of this period, \citet{thompson1967use} demonstrated that a consideration of extreme value distributions within a discriminal process leads to the Bradley--Terry model, and \citet{daniels1969round}, in a highly original paper, noted the links between the Bradley--Terry model and what might now be recognized as an undamped \PageRank\ \citep{page1999pagerank}. 

For further details of the development of the model up to this point \citet{david1963method} provides a thorough account of the paired comparison literature more generally, \citet{bradley1976science} and \citet{davidson1976bibliography} give interesting perspectives on the literature related to the Bradley--Terry model at the end of this period, and \citet{glickman2013introductory} is a highly readable account of the history, particularly as it pertains to the contribution of Zermelo. 

The next significant contributions to motivating the Bradley--Terry model came from \citet{henery1986interpretation} and \citet{joe1988majorization} in identifying the model as the result of maximizing an objective function subject to a suitable constraint. The later work \citep{joe1988majorization} seems to have been unaware of \citet{henery1986interpretation}, but provides a more complete presentation. As well as considering the Bradley--Terry model as a maximum entropy model and noting its relationship to an appropriate sufficient statistic, \citet{joe1988majorization} also explicitly notes the link to a maximum likelihood derivation. A number of motivations in this paper are based on game-style scenarios. Perhaps the most interesting paper related to this also comes from this period \citep{stern1990continuum}. In the context of the purpose of this work, \citet{mccullagh1993models} provides a particularly pertinent contribution at the end of this period, demonstrating how the Bradley--Terry model can be motivated from a geometric perspective, as well as how, under certain conditions, it is essentially equivalent to two other well-known models for permutations and from directional statistics respectively.

More recently \citet{slutzki2006scoring}, \citet{negahban2012iterative}, \citet{maystre2015fast} and \citet{selby2020citation} provide more detailed accounts of the link between the Bradley--Terry model and the limiting distribution of a Markov Chain, and thereby to an undamped \PageRank. The Social Choice literature provides an interesting perspective on this relationship, building on the approach of \citet{rubinstein1980ranking} to provide axiomatic justifications for ranking methods. \citet{slutzki2006scoring} is perhaps the most notable example in the present context.  

The paper proceeds by dividing the motivations up into six types: axiomatic; objective function maximization; discriminal processes; standard models; game scenarios; and quasi-symmetry and consistent estimators. These categorizations are somewhat arbitrary, and linkages exist across them which will be highlighted, but for the present purpose they provide a useful means to order the work. It begins with Section \ref{sec: axiomatic}, the discussion of axiomatic approaches, which takes as a starting point features that one might reasonably desire of a pairwise comparison model. A number are very closely linked and might even be thought of as restatements of the same idea, but the intuitions behind them differ sufficiently, as evidenced by their separate appearances in the literature, such that they are presented separately here. 

In Section \ref{sec: objective function maximisation}, the selection of a rating model is cast in the familiar framework of a constrained optimization, where an objective function is maximized or minimized subject to some plausible constraints. Section \ref{sec: discriminal processes} takes the context of Thurstone's discriminal processes, and discusses the distributions that lead to a Bradley--Terry model under this set-up, and how they might be motivated. In Section \ref{sec: standard models}, it is noted how the Bradley--Terry model is apparent in other well-known statistical models, as a conditional form of Rasch, Mallows $\phi$, von Mises--Fisher, hazard and network models. In Section \ref{sec: games}, some examples are introduced that derive from realistic game scenarios picking up on the highly intuitive nature of the model. In Section \ref{sec: quasi-symmetry}, the quasi-symmetry model is discussed, and is used to show how the often intuitive approaches that underlie a number of other popular rating methods can be related to Bradley--Terry and produce consistent estimators for the Bradley--Terry strength parameters. This also leads to noting the link to Barker's algorithm, a popular Markov chain Monte Carlo method. 

In each subsection, the reference given in the title is that of the earliest work linking the approach explicitly to the Bradley--Terry model, and the subsections are ordered chronologically by these. Where no reference is given, the link is believed to be novel. The sections are ordered with statistical interest and chronology in mind. 

In Section \ref{sec: discussion}, the natural questions of how these motivations are linked and the usefulness of motivating the model from diverse perspectives is addressed. The linkages are established with an examination of the Bradley--Terry model in the context of an exponential family of distributions. In demonstrating the usefulness of the approach, two illustrative examples are provided where it may be natural to use the model based on one motivation, but its application can be aided by considering it through another motivation.

Throughout the paper, $p_{ij}$ will be the probability of $i$ beating $j$ or for a preference for $i$ over $j$ given a comparison between $i$ and $j$, where $i,j \in T$ and $T$ is of size $n$. The $n \times n$ data matrix $C = [c_{ij}]$ will be the `competition' matrix of preferences or wins, such that $c_{ij}$ is the number of times $i$ was preferred over $j$. We define also $M=C + C^T$ as the symmetric matrix where $m_{ij}$ is the number of comparisons, or `matches' in British sports parlance, between $i$ and $j$. For the avoidance of doubt, no item is compared with itself, so that $c_{ii}=m_{ii}=0$ for all $i$. The observed number of wins for $i$ is denoted by $w_i=\sum_j c_{ij}$. Matrix $C$ is taken to be irreducible, that is, as described by \citet[p.29]{ford1957solution}: ``[I]n every possible partition of the objects into two non-empty subsets, some object in the second set has been preferred at least once to some object in the first set.'' This ensures that strength estimates are finite. It is not assumed that there are the same number of comparisons between any two items, nor indeed that the number of comparisons between any two items is non-zero. Shortened summation notation is used such that $\sum_{i,j}$ is taken to be $\sum^n_{i=1} \sum^n_{j=1}$ and $\sum_{i<j}$ is taken to be $ \sum^n_{j=1} \sum^{j-1}_{i=1}$. Where appropriate, the language of sports --- contests, scores, teams, wins --- is used to aid in providing clear interpretability, though the motivations may be analogized outside this context. 

For the convenience of readers, we end this Introduction with an index to the various Bradley--Terry model motivations that are included in this article: 

\begin{itemize}
	\item[\ref{sec: axiomatic}] Axiomatic motivations
	\begin{itemize}
		\item[\ref{transitivity of odds}] Transitivity of odds \citep{good1955marking}
		\item[\ref{sec: Luce}] Luce's Choice Axiom \citep{luce1959individual}
		\item[\ref{sec: reciprocity}] Reciprocity \citep{block1960contributions}
		\item[\ref{sec: sufficient statistic}] Wins as a sufficient statistic \citep{buhlmann1963pairwise}
	\end{itemize}
	\item[\ref{sec: objective function maximisation}] Objective function maximization
	\begin{itemize}
		\item[\ref{sec: maxent}] Maximum entropy with retrodictive criterion \citep{henery1986interpretation, joe1988majorization}
		\item[\ref{sec: geometric}] Geometric minimization \citep{mccullagh1993models}
		\item[\ref{sec: simplicity1}] Maximum definitional simplicity 1
		\item[\ref{sec: simplicity2}] Maximum definitional simplicity 2
	\end{itemize}
	\item[\ref{sec: discriminal processes}] Discriminal processes
	\begin{itemize}
		\item[\ref{sec: exponential}] Exponential distribution \citep[Holman and Marley as cited by][]{luce1965preference}
		\item[\ref{sec: extreme}] Extreme value distributions \citep{bradley1965another, thompson1967use}
	\end{itemize}
	\item[\ref{sec: standard models}] Standard models
	\begin{itemize}
		\item[\ref{sec: Rasch}] Rasch model \citep{andrich1978relationships}
		\item[\ref{sec: Mallows}] Mallows' $\phi$-model \citep{mccullagh1993models}
		\item[\ref{sec: von Mises--Fisher}] von Mises--Fisher distribution \citep{mccullagh1993models}
		\item[\ref{sec: CoxPH}] Cox proportional hazards model \citep{su2006connection}
		\item[\ref{sec: Network models}] Network models
	\end{itemize}
	\item[\ref{sec: games}] Game scenarios
	\begin{itemize}
		\item[\ref{sec: poisson scoring}] Poisson scoring \citep{audley1960stochastic, stern1990continuum}
		\item[\ref{sec: sudden death}] Sudden death \citep{stirzaker1999probability, vojnovic2015contest}
		\item[\ref{sec: accumulated win ratio}] Accumulated win ratio \citep{vojnovic2015contest}
		\item[\ref{sec: continuous time state transition}] Continuous time state transition
	\end{itemize}
	\item[\ref{sec: quasi-symmetry}] Quasi-symmetry and consistent estimators
	\begin{itemize}
		\item[\ref{sec: pagerank}] \PageRank\ \citep{daniels1969round}
		\item[\ref{sec: fair bets}] Fair bets \citep{daniels1969round}
		\item[\ref{sec: Wei--Kendall}] Wei--Kendall 
		\item[\ref{sec: RPI}] Ratings Percentage Index
		\item[\ref{sec: barker}] ``Winner stays on'' --- Barker's algorithm
	\end{itemize}
\end{itemize}

%% Next 3 lines are now superseded by the more concise index of sections 2-7 that appears above
%\renewcommand{\contentsname}{Motivations}
%\tableofcontents
%\addtocontents{toc}{\protect\setcounter{tocdepth}{2}}

\section{Axiomatic motivations} \label{sec: axiomatic}

It is sometimes possible to fix properties that we would desire of a model and use them to derive a unique model. In this section, we consider such properties that lead to the Bradley--Terry model.

\subsection{Transitivity of odds \citep{good1955marking}} \label{transitivity of odds}

Consider four teams $i,j,k,l$. Suppose that the probability that $j$ beats $k$ is greater than the probability that $j$ beats $l$,
\[p_{jk} > p_{jl},
\]
then it is intuitive to think that the probability that $i$ beats $k$ will be greater than the probability that $i$ beats $l$,
\[p_{ik} > p_{il}.
\]
Perhaps the simplest way to enforce this is by insisting on the transitivity of odds as \citet{good1955marking} proposes, that is
\[
\frac{p_{ij}}{p_{ji}} \times \frac{p_{jk}}{p_{kj}} = \frac{p_{ik}}{p_{ki}}.
\]
Alternatively, one might think of the same condition in the manner that \citet{luce1965preference} refers to it as the \emph{product rule}, where for any triple $(i,j,k)$ the probability of the intransitive cycle $i$ beats $j$, $j$ beats $k$, $k$ beats $i$ is the same as that of the intransitive cycle $i$ beats $k$, $k$ beats $j$, $j$ beats $i$, expressed
\[ 
p_{ij}p_{jk}p_{ki}=p_{ik}p_{kj}p_{ji} \quad \text{for all triplets }(i, j ,k).
\]
\citet{strang2022network} characterize this as an `arbitrage free' condition. Alternatively it may be understood as an energy conservation condition (on the log scale), and it is also known as Kolmogorov's criterion for reversibility of a Markov chain \citep{kolmogoroff1936theorie, kelly1979reversibility}. 

\citet{jech1983ranking} provides an alternative justification for the principle by considering estimating the odds of an item $i$ beating an item $k$ in the scenario where the comparison can only be made indirectly by comparing $i$ to $j$ and $j$ to $k$. If $i$ beats $j$ and $j$ beats $k$ then $i$ is taken to have beaten $k$. If $i$ loses to $j$ and $j$ loses to $k$ then $k$ is taken to have beaten $i$. For other result combinations ($i$ beats $j$ and $k$ beats $j$, or $j$ beats $i$ and $j$ beats $k$) judgment is reserved. In any given comparison, the probability that $i$ beats $k$ is thus $p_{ik}=p_{ij}p_{jk}$ and the probability that $k$ beats $i$ is thus $p_{ki}=p_{ji}p_{kj}$. Taking the ratio of these probabilities, the odds conform to the transitivity condition. \citet[p.246]{jech1983ranking} claims that this leads to the ``one and only one correct way of comparing the records of teams in an incomplete tournament", which seems a little bold, but the argument nevertheless demonstrates the intuitive appeal of the property.

Returning to how this criterion leads to the Bradley--Terry model, and following \citet{good1955marking}, it may alternatively be expressed as
\[
\log\frac{p_{ij}}{p_{ji}} + \log\frac{p_{jk}}{p_{kj}} = \log\frac{p_{ik}}{p_{ki}}.
\]
Letting $p_{ij}/p_{ji}=\exp(\tau(\theta_i, \theta_j))$, where $\theta_i$ can be thought of as a parameter summarizing the strength of $i$, then 
\[
\tau(\theta_i, \theta_j) + \tau(\theta_j,\theta_k) = \tau(\theta_i,\theta_k).
\]
Setting $\theta_j = \theta_i$, it may be noted that $\tau(\theta_i, \theta_i) = 0$ for all $i$. By setting $\theta_k=\theta_i$ it may be noted that $\tau$ is an antisymmetric function. Further, by differentiating with respect to $\theta_i$ it may be noted that the partial derivative of $\tau(\theta_i, \theta_j)$ with respect to $\theta_i$ is independent of $\theta_j$, so that $\tau(\theta_i, \theta_j)$ is some function of $\theta_i$ alone plus some function of $\theta_j$ alone, and since $\tau$ is antisymmetric it must be of the form 
\[
\tau(\theta_i,\theta_j) = t(\theta_i)-t(\theta_j). 
\]
Since $\theta_i$ is the strength of $i$ then $\tau(\theta_i,\theta_j)$ must be a monotone increasing function of $\theta_i$ and so $t(\theta_i)$ is a monotone increasing function of $\theta_i$ also. Therefore, $\lambda_i = t(\theta_i)$ is also a strength parameter for $i$ and 
\[
\frac{p_{ij}}{p_{ji}} = \exp(\lambda_i - \lambda_j) \quad \text{for all } i,j,
\]
giving the Bradley--Terry model.

\subsection{Luce's Choice Axiom \citep{luce1959individual}} \label{sec: Luce}

Let $p_S(i)$ be the probability that item $i$ is chosen from a set $S \subseteq T$, then a complete system of choice probabilities satisfies Luce's Choice Axiom if and only if for every $i$
\[
p_S(i)=\frac{p_{T}(i)}{\sum_{k \in S}p_{T}(k)} \quad .
\]

\citet{luce1959individual} introduces the axiom with the assertion that many choice situations are characterized by a multistage process, whereby a subset of the total choice set is selected, from which further subsets are selected iteratively, until a single choice is made from one of these subsets. \citet{luce1959individual} notes that for complex choices and a multistage process, the final result is likely to depend on these intermediate categorizations. However, for a simple decision and a two stage process, it is argued that the two-stage choice, reflected by the product $p_S(i)\sum_{k \in S}p_T(k)$, does not depend on $S$. By setting $S=T$, it is apparent that this product must be equal to $p_T(i)$. The Choice Axiom itself has been motivated by appealing to the decomposition of a full ranking model \citep{block1960contributions}, to invariance under uniform expansion of the choice set \citep{yellot1977relationship}, and, under specific assumptions, in a consideration of the utility of gambling \citep{luce2008utility}. 

A complete system satisfies the Choice Axiom if and only if there exist a set of numbers $\pi_1, \pi_2, \dots \pi_n$ such that for every $i$ and every set $S \subseteq T$
\[
p_S(i)=\frac{\pi_i}{\sum_{k \in S}\pi_k} \quad .
\]
In order to see this, let
\[
\pi_i = \kappa p_T(i), \quad \kappa>0,
\]
then by Luce's Choice Axiom
\begin{align*}
p_S(i) &= \frac{p_T(i)}{\sum_{k \in S}p_T(k)} \\
&= \frac{\kappa p_T(i)}{\sum_{k \in S}\kappa p_T(k)} \\
&= \frac{\pi_i}{\sum_{k \in S}\pi_k}.
\end{align*}
$\pi_i$ is unique up to a multiplicative constant since suppose there is another $\pi'_i$ satisfying this condition, then 
\[
\pi_i = \kappa p_T(i) = \frac{\kappa \pi'_i}{\sum_{k \in T}\pi'_k},
\]
and setting $\kappa'= \kappa / \sum_{k \in T}\pi'_k$ then $\pi_i = \kappa'\pi'_i$.

Taking $S$ to be the two-member set $\{i,j\}$ gives the Bradley--Terry model.

\subsection{Reciprocity \citep{block1960contributions}} \label{sec: reciprocity}
What might be thought of as an alternative expression of the Choice Axiom is noted in \citet{block1960contributions}. The idea is that the odds of $i$ beating $j$ should be equivalent to the ratio of strength parameters of $i$ and $j$,
\[
\frac{p_{ij}}{p_{ji}}=\frac{\pi_i}{\pi_j} \quad \text{for all }i,j \quad .
\]
Of course this condition can be framed in other familiar equivalent terms, either as detailed balance, more typically expressed as
\[
p_{ij}\pi_j=p_{ji}\pi_i \quad \text{for all }i,j,
\]
or that the irreducible, positive recurrent, aperiodic Markov chain for which $P=[p_{ij}]$ is the transition matrix is reversible, which itself is the case if and only if the transitivity condition of Section \ref{transitivity of odds} holds \citep{kelly1979reversibility}. The condition leads immediately to
\[
p_{ij} = \frac{\pi_i}{\pi_i+\pi_j} \quad.
\]
The relationship to Markov chains is discussed further elsewhere in this work. An explicit motivation in the context of a discrete Markov chain is introduced in Section \ref{sec: continuous time state transition}, and the discussion of Section \ref{sec: quasi-symmetry} is also relevant, in particular with the link to Barker's algorithm, a prominent Markov chain Monte Carlo method, discussed in Section \ref{sec: barker}.

\subsection{Wins as a sufficient statistic \citep{buhlmann1963pairwise}} \label{sec: sufficient statistic}

Define a statistical model for pairwise comparison where the probability that $i$ beats $j$ is independent of other pairwise comparisons conditional on strengths $\pi_i$ and $\pi_j$. Suppose $w_i=\sum_j c_{ij}$ are the wins gained by team $i$ and that the wins vector $\boldsymbol{w}=(w_1,w_2, \dots,w_n)^T$ is a sufficient statistic for the strength vector $\boldsymbol{\pi}=(\pi_1,\pi_2, \dots,\pi_n)^T$. 

Consider the comparison matrix $C=[c_{ij}]$ with $c_{kl},c_{lm}$ and $c_{mk}$ non-zero, for the triplet $(k,l,m)$ where, without loss of generality, $k<l<m$. Now consider an alternative $C'$ with $c'_{kl}=c_{kl}-1$, $c'_{lm}=c_{lm}-1$, $c'_{mk}=c_{mk}-1$ and $c'_{lk}=c_{lk}+1$, $c'_{ml}=c_{ml}+1$, $c'_{km}=c_{km}+1$, and all else the same. Then the wins vectors for the tournaments represented by $C$ and $C'$ are identical. If wins are a sufficient statistic for the strength parameters then the likelihood is dependent on $C$ only through $\boldsymbol{w}$, and so the likelihoods must also be identical. The likelihood is
\[
\prod_{i<j} \binom{m_{ij}}{c_{ij}} p_{ij}^{c_{ij}}(1-p_{ij})^{m_{ij}-c_{ij}},
\]
so that the log-likelihood, up to a constant term, is
\[
\sum_{i<j}c_{ij} \log \bigg ( \frac{p_{ij}}{1-p_{ij}} \bigg ) + m_{ij} \log (1-p_{ij}).
\]
Setting these equal for $C$ and $C'$, we get that
\begin{eqnarray*}
0 &=& (c_{kl}-c'_{kl})\log \frac{p_{kl}}{p_{lk}} \\
&& + (c_{lm}-c'_{lm})\log \frac{p_{lm}}{p_{ml}} \\
&& + (c_{mk}-c'_{mk})\log \frac{p_{mk}}{p_{km}},
\end{eqnarray*}
and so
\[
\log \frac{p_{kl}}{p_{lk}} + \log \frac{p_{lm}}{p_{ml}} + \log \frac{p_{mk}}{p_{km}} = 0,
\]
by the specifications of $c'_{kl},c'_{lm},c'_{mk}$, giving the Bradley--Terry model following the same argument as in Section \ref{transitivity of odds}.

\section{Objective function maximization} \label{sec: objective function maximisation}

It is a common procedure in quantitative analysis to identify an appropriate objective function and seek to maximize (or minimize) that function under certain plausible constraints. Indeed the familiarity of such procedures makes these motivations perhaps some of the most persuasive in the use of the Bradley--Terry model. 

\subsection{Maximum entropy with retrodictive criterion \citep{henery1986interpretation,joe1988majorization}} \label{sec: maxent}
In order to determine a functional form for the $p_{ij}$, we wish to select an appropriate objective function $S(p)$, a function of the probabilities $p_{ij}$, and then maximize this objective function subject to some appropriate constraint. 

The proposed constraint is that of the `retrodictive criterion', that the observed number of wins for each team is equal to the expected number of wins given the matches played. That is
\[
w_i = \sum_{j} c_{ij} = \sum_{j}m_{ij}p_{ij} \quad \quad \text{for all teams } i.
\]
A justification for this criterion was pithily expressed by \citet[p.280]{stob1984supplement} in summarizing the argument of \citet{jech1983ranking}: ``What sort of a claim is it that a team solely on the basis of the results should have expected to win more games than they did?'' This would seem to fail to appreciate the potential for bias from finite observations; nevertheless, it reflects the intuitive appeal of the condition.

Alternatively, the framework provided by \citet{Firth2022Alt3} offers a justification for the retrodictive criterion based on the equivalence of two intuitive formulations for rating in this setting. The first formulation proposes that given the pairwise win probabilities $p_{ij}$, an intuitive rating for a team $i$ would be the average win probability against all other competitors
\[
\Bar{p}_{i \cdot} = \frac{1}{n-1}\sum^n_{j=1, j \neq i}p_{ij}.
\]
The second formulation takes the ratio of observed wins for $i$ divided by the `effective matches' played by $i$, 
$w_i / m'_i$. Effective matches played, $m'_i$, is chosen to account for the strength of opposition. Any definition of $m'_i$ should meet two criteria. First, if the opponents played by $i$ have been strictly stronger (weaker) than average, then $m'_i$ is strictly less (greater) than $m_i$, the matches played by $i$, thus making the value of observed wins per effective matches played greater (less) than the value of observed wins per matches played. Second, observed wins per effective matches played, $w_i / m'_i$, is equal to the observed wins per matches played, $w_i / m_i$, in the case of a round-robin tournament, so that the rating accords with round-robin ranking. The simplest proposal meeting these two criteria is to scale each match played by the ratio of the probability of winning that match to the average probability of winning a match,
\[
m'_i = \sum_j m_{ij}p_{ij} / \Bar{p}_{i \cdot}.
\]
If we then set these two ratings, $\Bar{p}_{i \cdot}$ and $w_i / m'_i$ equal for all teams $i$, then we get the retrodictive criterion.

Turning to the objective function, the approach of maximizing entropy is common in statistical physics. Entropy is a measure of the uncertainty of a random variable. By maximizing it, roughly speaking, the assumptions in the model are minimized. \citet{jaynes1957information} influentially advocated for the choice of entropy in a broader range of statistical settings, building on the ideas from information theory of \citet{shannon1948mathematical}. \citet{good1963maximum} provides further discussion noting ``[t]he mere fact that the principle of maximum entropy generates classical statistical mechanics, as a null hypothesis, would be sufficient reason for examining its implications in mathematical statistics.'' \citet{luce1959individual}, on the other hand, casts doubt on its applicability to choice contexts. 

In this setting, the entropy is defined as
\begin{eqnarray*}
S(p) &=& -\sum_{i,j} m_{ij}p_{ij}\log p_{ij} \\
&=& -\sum_{i<j}m_{ij}(p_{ij} \log p_{ij} + (1-p_{ij}) \log (1-p_{ij})) .
\end{eqnarray*}
We maximize the entropy subject to the retrodictive criterion using the method of Lagrange multipliers,
\[
    \mathcal{L}(p,\boldsymbol{\eta}) = S(p) - \sum_{i=1}^n \eta_i\bigg(\sum_{j=1}^n (m_{ij}p_{ij} - c_{ij})\bigg),
\]
and setting $\frac{\partial \mathcal{L}}{\partial p_{ij}}=0$ for all $p_{ij}$ in the normal way gives that
\begin{align*}
\frac{\partial S(p)}{\partial p_{ij}} &=\frac{\partial}{\partial p_{ij}}\sum_{r=1}^n \eta_r \bigg(\sum_{s=1}^n (m_{rs}p_{rs} - c_{rs})\bigg) \quad \text{for all }i,j.
\end{align*}
So for all $i, j$ such that $m_{ij} \neq 0$,
\[
- \log p_{ij} + \log(1-p_{ij}) = \eta_i - \eta_j , 
\]
or equivalently
\[
p_{ij} = \frac{\pi_i}{\pi_i + \pi_j},
\]
where $\pi_i = \exp(-\eta_i)$, and it can readily be checked by differentiating $S(p)$ that this is a maximum. 

\subsection{Maximum likelihood estimation with retrodictive criterion \citep{joe1988majorization}} \label{sec: maxlikelihood}

Maintaining the retrodictive criterion of Section \ref{sec: maxent}, we might consider the likelihood as an alternative objective function to maximize. This is consistent with the use of likelihood-based information criteria, such as AIC and BIC, for model choice. Suppose the probability of $i$ being preferred to $j$ is given by
\[
p_{ij}=f(\lambda_i,\lambda_j),
\]
where $\lambda_i$ and $\lambda_j$ are real-valued parameters describing the strength of items $i$ and $j$, and $f :  \mathbb{R}\times\mathbb{R}\rightarrow [0,1]$. Then the likelihood function is given by
\begin{eqnarray*}
    L(\boldsymbol{\lambda}) &=& \prod_{i<j}\binom{m_{ij}}{c_{ij}}p_{ij}^{c_{ij}}(1-p_{ij})^{m_{ij}-c_{ij}} \\
    &=& \prod_{i<j}\binom{m_{ij}}{c_{ij}}p_{ij}^{c_{ij}}p_{ji}^{c_{ji}},
\end{eqnarray*}
and the log-likelihood function, ignoring the constant term, is
\[
l(\boldsymbol{\lambda}) 
= \sum_{i<j}{c_{ij}}\log(p_{ij}) + {c_{ji}}\log(p_{ji}).
\]
At an extreme point of the log-likelihood, for all $k$,
\[
0 = \frac{\partial}{\partial \lambda_k}l(\boldsymbol{\lambda}) \nonumber
 = \sum_{j} c_{kj}\frac{\partial }{\partial \lambda_k}\log(p_{kj}) + c_{jk}\frac{\partial }{\partial \lambda_k}\log(p_{jk}).
\]
Considering the constraint we note that
\begin{align*}
0 = \sum_{j} c_{kj} - m_{kj}p_{kj} \nonumber
&= \sum_{j} c_{kj} - (c_{kj}+c_{jk})p_{kj} \nonumber \\
&= \sum_{j} c_{kj}(1-p_{kj}) - c_{jk}p_{kj} \nonumber \\
&= \sum_{j} c_{kj}(1-p_{kj}) - c_{jk}(1-p_{jk}),
\end{align*}
and so there is an extreme point where
\begin{align*}
    \frac{\partial }{\partial \lambda_k}\log(p_{kj}) &= (1-p_{kj}) \quad \text{and}\\
    \frac{\partial }{\partial \lambda_k}\log(p_{jk}) &= -(1-p_{jk}),
\end{align*}
which gives
\begin{align*}
    \frac{\partial p_{kj}}{\partial \lambda_k}&= p_{kj}(1-p_{kj}) \quad \text{and}\\
    \frac{\partial p_{jk}}{\partial \lambda_k} &= -p_{jk}(1-p_{jk}).
\end{align*}
Solving these separable differential equations for $p_{ij}$ gives
\begin{align*}
p_{ij}&=\frac{e^{(\lambda_i-\lambda_j)}}{1+e^{(\lambda_i-\lambda_j)}} \\
& = \frac{\pi_i}{\pi_i + \pi_j}
\end{align*}
where $\pi_i = e^{\lambda_i}$, and, as before, this is a maximum since the log-likelihood is strictly concave. So that the Bradley--Terry model is the likelihood-maximizing model.

\subsection{Geometric minimization \citep{mccullagh1993models}}\label{sec: geometric}

If one were to conceive of the rating of $n$ items under a geometric interpretation, a natural general framing might start by representing the observed results as vectors in some $n$-space. A rating vector can then be taken as the vector that minimizes some aggregate quantity with respect to these observed result vectors, where the quantity is smaller when the rating vector is more concordant with the results. \citet{mccullagh1993models} presents just such a framing with the outcome and rating vectors confined to a $n$-sphere, taken to be of unit radius for convenience. For example, in a five-team tournament consisting of competitors $A, B, C, D, E$ then a win for $D$ over $B$ would be represented by the result vector $\mathbf{x} = (0, -1/\sqrt{2}, 0, 1/\sqrt{2}, 0)$. 

With both the rating vector and observed result vectors lying on the unit sphere, a natural quantity to seek to minimize is the angle between the rating vector and an observed result vector, or equivalently maximizing the cosine of the angle as expressed through the dot product of the vectors, $\mathbf{x}\cdot\mathbf{\lambda}$. Note that this is equivalent to minimizing the squared Euclidean distance between the points on the sphere since
\[
 \mid \mid \mathbf{x} - \mathbf{\lambda} \mid \mid^2 = 2 - 2\mathbf{x}\cdot\mathbf{\lambda}.
\]
So to find our rating vector $\mathbf{\lambda}$, we would sum the dot product over all observed results and select $\mathbf{\lambda}$ such that it maximizes this quantity. 

In the notation used in this paper, and keeping the unit radius, any result vector $\mathbf{x_{ij}}$ representing a win for $i$ over $j$ will have value $1/\sqrt{2}$ in the $i$th position, $-1/\sqrt{2}$ in the $j$th position and zero elsewhere. The sum over all such results is therefore
\[
\frac{1}{\sqrt{2}}\sum_{i,j} c_{ij}\left( \lambda_i - \lambda_j \right) = \frac{1}{\sqrt{2}}\sum_{i,j}\lambda_i(c_{ij}-c_{ji}),
\]
which is the form of the likelihood maximization that gives the Bradley--Terry rating (see Section \ref{sec: exp family} for further details). Thus, a geometric interpretation of rating where one minimizes the aggregate angles between result vectors and a rating vector on a sphere returns the Bradley--Terry ratings. One nice feature of this motivation is the ready extendability to scenarios of differing numbers of competitors in each contest, while maintaining consistency with Bradley--Terry in the pairwise contest case. This is discussed further in Section \ref{sec: sports ranking}.

\citet{mccullagh1993models} also demonstrates the link to Mallows' $\phi$-model and the von Mises--Fisher distribution. These are presented in Section \ref{sec: standard models}.

\subsection{Maximum definitional simplicity 1}\label{sec: simplicity1}

Often when selecting a model, transparency and interpretability are desirable features. This may be especially so in contexts where fairness of a ranking system are a consideration. These sort of contexts are common in pairwise comparison with the methods being used to perform activities like ranking sports teams or players \citep{zermelo1929berechnung, Firth2022Alt3} or in educational assessment \citep{pollitt2012method}. Therefore, there may be a legitimate desire for definitionally simpler, more intuitive models. It is thus appealing to consider how one might select a model with the goal of maximizing definitional simplicity.

Suppose one wished to determine a ranking by defining a probability for the preference for $i$ over $j$ related only to positive real-valued strength parameters $\pi_i$ and $\pi_j$ respectively, 
\[
p_{ij}=f(\pi_i,\pi_j).
\]

A reasonable set of criteria for this function would be:
\begin{enumerate}
\item $f: \mathbb{R}^{+}\times\mathbb{R}^{+} \rightarrow [0,1]$,
\item $f(\pi_i,\pi_j) = \frac{1}{2}$ when $\pi_i=\pi_j$,
\item $\lim_{\pi_i\to 0, \pi_j \text{ fixed}} f(\pi_i,\pi_j)=0$,
\item $\lim_{\pi_j\to 0, \pi_i \text{ fixed}} f(\pi_i,\pi_j)=1$,
\item $\lim_{\pi_i\to\infty, \pi_j \text{ fixed}} f(\pi_i,\pi_j)=1$,
\item $\lim_{\pi_j\to\infty, \pi_i \text{ fixed}} f(\pi_i,\pi_j)=0$.
\end{enumerate}
where  $\mathbb{R}^+$ is taken to be the set of positive real numbers not including zero.

Assume that the simplest set of functions are those that may be defined solely using the four basic operators ($+, -, \times, \div$), and that any measure of the simplicity of a function is a strictly decreasing function of the number of these operators used. In this setup, maximizing definitional simplicity of a function is equivalent to minimizing the number of basic operators in its definition. Bracketing anywhere, used in the conventional sense, to identify a functional sub-clause, is allowed without increasing or reducing simplicity. Constants are also allowed in place of parameters without increasing or reducing simplicity. In the language of Computer Science, this is equivalent to defining simplicity by the minimum number of floating point operations (flops).

First, we note that no function that is equivalent to a constant can meet all criteria simultaneously. The only functions with zero or one operator that meet criterion 5 are constants $f(\pi_i,\pi_j)=1$ or equivalents (for example, $f(\pi_i,\pi_j)=\pi_i / \pi_i$). Therefore, at least two operators must be required. Criterion 5 implies that the operator $\div$ is employed as otherwise the function would be a constant or the limit would be infinite in absolute value, violating the other criteria . So, if there is a solution with exactly two operators, then it must be of the form $f(\pi_i,\pi_j)=g(\pi_i,\pi_j) \div h(\pi_i,\pi_j)$. In order to ensure that there are only two operators in total, either $g$ or $h$ must have no operator, and therefore be equal to one of the parameters or to a constant. The other must be a single-operator function involving $+$ or $-$ in order to meet criterion 5 without being equivalent to a constant. From criterion 3, it must be that $g(\pi_i,\pi_j)=\pi_i$ and then from criterion 5, $h$ must take $\pi_i$ as one of its terms. Criterion 6 implies that the other term in $h$ is $\pi_j$ and criterion 2 then implies that $h(\pi_i,\pi_j)=\pi_i+\pi_j$. This gives $f(\pi_i,\pi_j)=\pi_i \div (\pi_i+\pi_j)$, which meets all the required criteria. It may be noted that not all the criteria were required for its unique derivation, and that other subsets of the criteria may be used to derive the same result. That is to say that 
\[
 p_{ij} = \frac{\pi_i}{\pi_i + \pi_j}
\]
will be the unique simplicity maximizer under a number of different subsets of the plausible criteria.

\subsection{Maximum definitional simplicity 2}\label{sec: simplicity2}

Given positive-valued strength parameters $\pi_i$ and $\pi_j$ for $i$ and $j$ respectively, one may want to consider a model where the probability of $i$ being preferred to $j$ is a function of the ratio $x_{ij}=\pi_i / \pi_j$,
\[
p_{ij}=f(x_{ij}).
\]

A reasonable set of criteria for this function would then be:
\begin{enumerate}
\item $f: \mathbb{R}^{+} \rightarrow [0,1]$,
\item $f(1) = \frac{1}{2}$,
\item $\lim_{x\to 0} f(x)=0$,
\item $\lim_{x\to\infty} f(x)=1$,
\end{enumerate}

Proceeding in a similar fashion to the previous section, constant functions cannot meet all criteria. The only function including exactly zero or one flop that meets criterion 4 is a constant. Considering a function with two operators and again considering criterion 4, then it must be that the operator $\div$ is employed as otherwise the function would be a constant or the limit would be infinite in absolute value. So if there is a solution with exactly two operators then it must be of the form $f(x)=g(x) \div h(x)$. In order to ensure that there are only two operators in total, either $g$ or $h$ must have no operator, and therefore be equal to one of the parameters or to a constant. The other must be a single-operator function involving $+$ or $-$ in order to meet criterion 4 without being equal to a constant. Criterion 3 implies that $g(x)=x$, and criterion 2 then tells us that $h(x)=1+x$. Thus
\[
f(x)=\frac{x}{1+x},
\]
giving
\[
 p_{ij} = \frac{\pi_i}{\pi_i + \pi_j}.
\]

\section{Discriminal processes} \label{sec: discriminal processes}

Consider a scenario where the strength of each of two entities in a given pairwise interaction is independently observed with error and then compared. The item with the greater observed strength is preferred. This is the model of Thurstone's `discriminal processes' \citep{thurstone1927law}. Denote the observed strength of $i$ as $b_i$ with `true' strength $\lambda_i$, so that $b_i = \lambda_i + \epsilon_i$, where $\epsilon_i$ is an error term. Item $i$ is preferred to item $j$ if and only if $b_i > b_j$. Taking the error to be Gaussian, as Thurstone himself did, leads to what is commonly known as the Thurstone--Mosteller model \citep{thurstone1927law, mosteller1951remarks}, but the set up may also be used to motivate the Bradley--Terry model by considering alternative distributions for $b_i$.

\subsection{Exponential distribution \citep[Holman and Marley as cited by][p.338] {luce1965preference}} \label{sec: exponential}

Suppose ${b_i}$ and $b_j$ follow independent exponential distributions whose expected values are given by $\pi_i$ and $\pi_j$ respectively with the cdf,
\[
F_i(x) = 1- \mathrm{e}^{-\frac{x}{\pi_i}}, \quad x \in \mathbb{R}^+.
\]
Then, with $F'$ denoting the pdf, the probability that $i$ is preferred to $j$ in a pairwise comparison is
\begin{align*}
p_{ij}  &= \int_{0}^\infty F_j(x)                        F_i'(x)\mathrm{d}x \\
        &= \int_{0}^\infty \left(1- \mathrm{e}^\frac{x}{\pi_j}\right)\frac{1}{\pi_i}\mathrm{e}^{-\frac{x}{\pi_i}} \mathrm{d}x \\
        &= 1 - \frac{1}{\pi_i\left(\frac{1}{\pi_i} + \frac{1}{\pi_j}\right)}\int_{0}^\infty \left(\frac{1}{\pi_i} + \frac{1}{\pi_j}\right) \mathrm{e}^{-\left(\frac{1}{\pi_i} + \frac{1}{\pi_j}\right)x} \mathrm{d}x \\
        &= 1 - \frac{\pi_j}{\pi_i + \pi_j} \\
        &= \frac{\pi_i}{\pi_i + \pi_j} \quad .
\end{align*}

\subsection{Extreme value distributions \citep{bradley1965another, thompson1967use}} \label{sec: extreme}

\citet{thompson1967use} provide a rationale for a broader class of distributions that lead to a Bradley--Terry model under a discriminal process. Based on ideas from Psychology, sensations are hypothesized to be a result of a large number of stimuli. These stimuli are modeled as having independent identical distributions $G(x)$. One might then consider the distribution of the resultant sensation. 

Two intuitive possibilities would be to model the distribution of the sensation $F(x)$ either as the average of those stimuli or the maximum of those stimuli. Taking the average gives a normal distribution for $F(x)$ in the limit, leading to a Thurstone--Mosteller comparison model. Taking the maximum of the stimuli, in the limit, gives, by extreme value theorem \citep{fisher1928limiting,gnedenko1943distribution, gumbel1958statistics}, one of three distributions for $F(x)$ --- Gumbel, Weibull, or Fr\'echet --- depending on the underlying stimuli distribution $G(x)$, leading to a Bradley--Terry comparison model. The Gumbel is the most notable of these, being the sensation distribution for stimuli distributions such as the normal, lognormal, logistic, and exponential. 

While \citet{thompson1967use} provided a clear motivation for considering such models and do not assume that the underlying stimuli distributions need have the same location parameters for $i$ and $j$, \citet{lehmann1953power} had previously considered a family of distributions in the context of the power of rank tests of the form $F_{X_i}(x;\pi_i) = G^{\pi_i}(x)$, where $G(x)$ is itself a distribution function. \citet{bradley1965another} discussed this family of distributions with respect to the Bradley--Terry model. As \citet{bradley1976science} notes, if $G(x)$ is a distribution function, and $X_i$ is the random variable relating to a sensation $i$, with distribution function 
\[
\mathbb{P}(X_i \leqslant x) = G^{\pi_i}(x) ,
\]
where $\pi_i>0$, then comparing sensations $i$ and $j$, 
\begin{align*}
p_{ij} &= \mathbb{P}(X_i > X_j) \\
&= \int_{x_i>x_j} dG^{\pi_i}(x_i)dG^{\pi_j}(x_j) \\
&= \frac{\pi_i}{\pi_i + \pi_j}.
\end{align*}

\subsubsection{Gumbel distribution \citep{thompson1967use}} 

Suppose ${b_i}$ follows a Gumbel distribution with mean $\lambda_i$. Then
\[
\boldsymbol{\text{Pr}}(b_i \leq x) = F_i(x) = \exp{(-\pi_i \mathrm{e}^{-\alpha x})} ,
\]
for $~ x \in \mathbb{R}$ and parameter $\alpha > 0$, where $\pi_i = \mathrm{e}^{\alpha \lambda_i - \gamma}$, with $\gamma$ the Euler--Mascheroni constant. Define $F_{i+j}(x) = \exp{(-(\pi_i+\pi_j) \mathrm{e}^{-\alpha x})}$ then the probability that $i$ is preferred to $j$ in a pairwise comparison is
\begin{align*}
p_{ij} &= \int_{-\infty}^\infty F_j(x)F_i'(x)\mathrm{d}x \\
        &= \int_{-\infty}^\infty \exp{(-\pi_j\mathrm{e}^{-\alpha x})}\alpha \pi_i \exp{(-\alpha x - \pi_i\mathrm{e}^{-\alpha x})} \mathrm{d}x \\
        &= \frac{\pi_i}{\pi_i + \pi_j}\int_{-\infty}^\infty F_{i+j}'(x) \mathrm{d}x \\
        &= \frac{\pi_i}{\pi_i + \pi_j} \quad .
\end{align*}

\subsubsection{Weibull distribution \citep{thompson1967use}}

Suppose ${b_i}$ follows a Weibull distribution
\[
\mathbb{P}(b_i \leq x) = F_i(x) = 1 - \exp{(-(x/\lambda_i)^\alpha)},
\]
for $ x \in \mathbb{R}^+$ and parameter $\alpha > 0$. Then the probability that $i$ is preferred to $j$ in a pairwise comparison is
\begin{align*}
p_{ij} &= \int_{0}^\infty F_j(x) F_i'(x)\mathrm{d}x \\
        &= \int_{0}^\infty [1-\exp{(-(x/\lambda_j)^\alpha)}]\frac{\alpha}{\lambda_i} \\
        & \quad \quad \quad \quad (x/\lambda_i)^{\alpha-1}\exp{(-(x/\lambda_i)^{\alpha})} \mathrm{d}x \\
        &= 1 - \int_{0}^\infty \frac{\alpha}{\lambda_i}(x/\lambda_i)^{\alpha-1} \\
        & \quad \quad \quad \quad \quad \quad \quad \quad \exp{(-(x/\lambda_j)^{\alpha})-(x/\lambda_i)^\alpha)} \mathrm{d}x \\
        &= 1 - \frac{\lambda_j^\alpha}{\lambda_i^\alpha + \lambda_j^\alpha} \int_{0}^\infty \frac{\alpha}{\lambda_i\lambda_j}(x/\lambda_i\lambda_j)^{\alpha-1}(\lambda_i^\alpha + \lambda_j^\alpha) \\
        & \quad \quad \quad \quad \quad \quad \quad \quad \exp{(-(x/\lambda_i\lambda_j)^{\alpha}(\lambda_i^\alpha + \lambda_j^\alpha))} \mathrm{d}x \\
        &= \frac{\pi_i}{\pi_i + \pi_j}\quad ,
\end{align*}
where $\pi_i=\lambda_i^\alpha$.

\subsubsection{Fr\'echet distribution \citep{thompson1967use}}

Suppose ${b_i}$ follows a Fr\'echet distribution
\[
\mathbb{P}(b_i \leq x) = F_i(x) = \exp{(-\pi_ix^{-\alpha})} 
\]
for $x \in \mathbb{R}^+$ and parameter $\alpha > 0$. Then the probability that $i$ is preferred to $j$ in a pairwise comparison is
\begin{align*}
p_{ij} &= \int_{0}^\infty F_j(x)                         F_i'(x)\mathrm{d}x \\
        &= \int_{0}^\infty \exp{(-\pi_jx^{-\alpha})} \frac{\pi_i \alpha}{x^{\alpha+1}}\exp{(-\pi_ix^{-\alpha})} \mathrm{d}x \\
        &= \frac{\pi_i}{\pi_i + \pi_j}\int_{0}^\infty \alpha \frac{\pi_i + \pi_j}{x^{\alpha+1}} \exp{(-(\pi_i+\pi_j)x^{-\alpha})} \mathrm{d}x \\
        &= \frac{\pi_i}{\pi_i + \pi_j}\quad .
\end{align*}

\section{Standard models} \label{sec: standard models}

A number of models familiar to statisticians may be related to the Bradley--Terry model by considering conditional forms. In Section \ref{sec: geometric}, we noted how \citet{mccullagh1993models} demonstrated links to Mallows' $\phi$-model and the von Mises--Fisher distribution. We expand on those links here and also discuss the relation to three more models familiar to statisticians.

\subsection{Rasch model \citep{andrich1978relationships}} \label{sec: Rasch}

Let $X_{vi}$ be a binary random variable, representing the outcome of a test $v$ taken by candidate $i$, where $X_{vi}=1$ represents passing the test, and $X_{vi}=0$ denotes failure. Under the Rasch simple logistic model \citep{rasch1960probabilistic, rasch1961general} the probability of the outcome $X_{vi}=1$ is taken to be
\[
\mathbb{P}(X_{vi}=1)=\frac{e^{\lambda_i-\delta_v}}{1+e^{\lambda_i-\delta_v}},
\]
where $\lambda_i$ represents the ability of candidate $i$ and $\delta_v$ the difficulty of test $v$. 

There are two conceptualizations by which we might derive the Bradley--Terry model from this. First, as \citet{andrich1978relationships} notes, if we take $p_{ij}$ to be 
\[
\mathbb{P}(i \text{ passes a test } v \mid \text{exactly one of }i \text{ and }j\text{ pass the test } v),
\]
then since
\[
\mathbb{P}(X_{vi}=1,X_{vj}=0)=\frac{e^{\lambda_i-\delta_v}}{(1+e^{\lambda_i-\delta_v})(1+e^{\lambda_j-\delta_v})},
\]
and
\[
\mathbb{P}(X_{vi}+X_{vj}=1)=\frac{e^{\lambda_i-\delta_v}+e^{\lambda_j-\delta_v}}{(1+e^{\lambda_i-\delta_v})(1+e^{\lambda_j-\delta_v})}
\]
then conditional on being able to discern that one of the test-takers has performed better based on the binary test outcome and taking their test outcomes to be independent conditional on their abilities and the test difficulty then the probability that $i$ has beaten $j$ is
\[
p_{ij}=\frac{\mathbb{P}(X_{vi}=1,X_{vj}=0)}{\mathbb{P}(X_{vi}+X_{vj}=1)}=\frac{e^{\lambda_i}}{e^{\lambda_i}+e^{\lambda_j}}  = \frac{\pi_i}{\pi_i + \pi_j},
\]
where $\pi_i = e^{\lambda_i}$. 

Second, we might more directly consider that in comparing $i$ with $j$ we are setting a test for $i$ of difficulty equal to the strength of the comparator $\lambda_j$ (or equivalently setting a test for $j$ of difficulty equal to the strength of the comparator $\lambda_i$), so that
\[
p_{ij}=\frac{e^{\lambda_i - \lambda_j}}{1+e^{\lambda_i - \lambda_j}}  = \frac{\pi_i}{\pi_i + \pi_j}.
\]

\subsection{Mallows' $\phi$-model \citep{mccullagh1993models}}\label{sec: Mallows}

\citet{mallows1957non} discusses models on the space of permutations. In the context of this paper, a permutation might equivalently be thought of as a ranking. The simplest of these model families is Mallows' $\phi$-model,
\[
p(\mathbf{x})=K_{\phi}\exp{\{-\phi d(\mathbf{x} , \lambda)\}},
\]
where $d$ is a distance measure between an observed permutation $\mathbf{x}$ and the `modal permutation' $\lambda$, $\phi$ is a concentration parameter and $K_{\phi}$ is a constant of proportionality. Thus, in maximizing the likelihood of the model given observed permutations, the modal permutation is the permutation that has the minimum aggregate distance to the observed permutations. In considering distances on permutations or ranks, the Spearman rank correlation coefficient is a natural candidate and is the one considered here.

\citet{mccullagh1993models} notes that the use of ordinal numbers to represent ranks, while a strong norm, is somewhat arbitrary. He proposes that the ranks are transformed such that 
\[
k' = \frac{k - (n+1)/2}{\sqrt{n(n^2-1)/12}},
\]
where $k$ is an integer from $1$ to $n$ representing a rank. With this transformation, the rank permutations are projected onto the unit sphere. For this paper, we consider a `modal rating' rather than a modal permutation or ranking and thus take the negative of this transformation to ensure that higher-ranked items have higher value. For example, a rank vector $(2,3,1,4)$, expressing that item 1 came second, item 2 third etc. would be transformed to the vector $\frac{1}{2\sqrt{5}
}(1, -1, 3, -3)$. The observed pairwise results may be projected onto the unit sphere with the pairwise ranking of a win for $i$ over $j$ represented by a vector with the value $1/\sqrt{2}$ in the $i$th position, $-1/\sqrt{2}$ in the $j$th position and zero elsewhere. 

If the pairwise results are represented in this way and the distance measure is taken to be the Spearman rank correlation coefficient, which is equivalent to the squared Euclidean distance, then following the argument from Section \ref{sec: geometric}, the Mallows' $\phi$-model becomes equivalent in form to the Bradley--Terry model, with the modal permutation vector equal to the vector of Bradley--Terry ratings. \citet{mccullagh1993models} notes that, strictly speaking, the models are not equivalent. Under the proposed ranking transformation, the Mallows' $\phi$-model is defined on the sample space of permutations represented on the unit sphere and has a ranking as the parameter $\lambda$, whereas the Bradley--Terry model is defined on the sample space of pairwise unit vectors and takes $\lambda$ to be any point on the unit sphere.

\subsection{von Mises--Fisher distribution \citep{mccullagh1993models}}\label{sec: von Mises--Fisher}
The von Mises--Fisher distribution \citep{von1918uber, fisher1953dispersion} is a well-known model in directional statistics. It defines a probability density for a random $n$-dimensional unit vector $\mathbf{x}$ as 
\[
p(\mathbf{x};\mathbf{\lambda})= C_{\kappa}\exp{\{ \kappa \mathbf{x} \cdot \mathbf{\lambda} \} }.
\]
As discussed in Sections \ref{sec: geometric} and \ref{sec: Mallows}, if we take the pairwise result outcomes and the Bradley--Terry rating vector to be defined on the unit sphere then this takes the same form as the Bradley--Terry model. As with the Mallows' model, \citet{mccullagh1993models} notes that, strictly speaking, they are not equivalent due to being defined on different sample spaces. In this case, the von Mises--Fisher distribution is defined on the continuous sample space of the unit sphere, whereas the Bradley--Terry model takes the pairwise rankings projected onto the unit sphere as its sample space.

\subsection{Cox proportional hazards model \citep{su2006connection}} \label{sec: CoxPH}

Consider a proportional hazards model \citep{cox1972regression} on random variables $T_i$ with hazard function given by
\[
h_i(t) = h(t)\pi_i.
\]
Thus the hazard rate for object $i$ is given by a multiplicative factor $\pi_i$. Then 
\begin{align*}
    &\mathbb{P}(T_i < T_j) \\
    & \quad = \int_0^{\infty} F_{T_i}(t)f_{T_j}(t) \,dt \\
    &\quad = \int_0^{\infty} \left( 1 - \exp \left\{-\int_0^t h(x)\pi_i \,dx \right\} \right) \\
    & \quad \quad \quad \quad \quad \quad h(t)\pi_j \exp\left\{-\int_0^t h(x)\pi_j \,dx\right\} \,dt \\
    &\quad = 1 - \int_0^{\infty} h(t)\pi_j \exp\left\{-(\pi_i + \pi_j)\int_0^t h(x) \,dx\right\} \,dt  \\
    &\quad = 1 - \frac{\pi_j}{\pi_i + \pi_j} \\
    &\quad = \frac{\pi_i}{\pi_i + \pi_j}.
\end{align*}
Further, as \citet{su2006connection} note, if a stratified proportional hazards model is used such that each stratum represents a different match with
\[
h_{ij}(t) = h_{s_{ij}}(t) \pi_i,
\]
where $s_{ij}$ is the stratum for a match between $i$ and $j$ then the contribution to the partial likelihood from the random variables $T_i$ and $T_j$ with the event $\{T_i < T_j\}$ is $\pi_i / (\pi_i + \pi_j)$.

\subsection[Network models]{Network models} \label{sec: Network models}
Consider a binary directed network $Y$, with an edge $i \rightarrow j$ taking the value $y_{ij}$. A common class of models in network analysis takes a conditional independence approach, assuming that the value of any directed edge is independent of all other edge values given an appropriate set of parameters. In a generalized form for the current purposes it can be expressed as
\begin{align*}
\mu_{ij} &= \mathbb{P}(y_{ij}=1 ) \\
\text{logit}(\mu_{ij} ; \delta_i, \gamma_j, f_{ij}) &= \delta_i + \gamma_j + f_{ij},
\end{align*}
where $\delta_i$ and $\gamma_j$, sometimes referred to as \emph{sociality} and \emph{attractivity} parameters \citep{krivitsky2009representing}, reflect the heterogeneity of out-degree and in-degree respectively, and $f_{ij}=f(i,j)$ is a symmetric function capturing the propensity for an edge in either direction to exist. For example, \citet{hoff2002latent} take $f(i,j)$ to be the Euclidean distance between points associated with $i$ and $j$ in a latent space but note that $f(i, j)$ could be any distance measure satisfying the triangle inequality $f(i,j) \leqslant f(i,k) + f(k,j)$. Often models also incorporate a term of the form $\beta^T {x_{ij}}$ within $f(i,j)$ , where ${x_{ij}}$ is a vector of pair-specific characteristics, in order to capture known homophilies.

Applying the conditional independence assumption and looking at the probability of an edge being present only in the direction $i \rightarrow j$ and not the $j \rightarrow i$ direction,
\begin{align*}
&\mathbb{P}(y_{ij}=1, y_{ji}=0 ; \delta_i, \delta_j, \gamma_i, \gamma_j, f_{ij}) \\
& \quad = \frac{ e^{\delta_i + \gamma_j + f_{ij}}}{(1 +e^{\delta_i + \gamma_j + f_{ij}})(1+e^{\delta_j + \gamma_i + f_{ij}})} ,
\end{align*}
so
\begin{align*}
&\mathbb{P}(y_{ij}=1 \mid y_{ij}+y_{ji}=1; \delta_i, \delta_j, \gamma_i, \gamma_j, f_{ij}) \\
& \quad = \frac{ e^{\delta_i + \gamma_j + f_{ij}}}{e^{\delta_i + \gamma_j + f_{ij}} + e^{\delta_j + \gamma_i + f_{ij}}} \\
& \quad = \frac{e^{\delta_i - \gamma_i}}{e^{\delta_i - \gamma_i} + e^{\delta_j - \gamma_j}} \\
& \quad = \frac{e^{\lambda_i}}{e^{\lambda_i} + e^{\lambda_j}} \\
& \quad = \frac{\pi_i}{\pi_i + \pi_j},
\end{align*}
where $\pi_i=e^{\lambda_i}$ and $\lambda_i = \delta_i - \gamma_i$. If $Y$ is considered as a tournament matrix with a directed edge $i \rightarrow j$ indicating $i$ beats $j$, then \emph{sociality} is a team's propensity for winning and  \emph{attractivity} the propensity for losing so that assessing the strength of a team as the difference between these is readily intuitive.

\section{Game scenarios} \label{sec: games}

The Bradley--Terry model has frequently been associated with an analysis of sport. So it is perhaps not surprising that there are a number of game scenarios in which the model may be very naturally motivated. Some of these are presented here.

\subsection{Poisson scoring \citep{audley1960stochastic, stern1990continuum}} \label{sec: poisson scoring}

Consider two teams $i$ and $j$ who score according to independent Poisson processes $X_i(t)$ and $X_j(t)$ with rate parameters $\pi_i$ and $\pi_j$ respectively. The winner is the first team to score. Then by Poisson thinning, for any value of $t$,
\begin{equation*}
    \begin{split}
        p_{ij} = \mathbb{P}(X_i(t)=1 \mid X_i(t)+X_j(t)=1) = \frac{\pi_i}{\pi_i + \pi_j} \quad .
    \end{split}
\end{equation*}

\citet{audley1960stochastic} presents an argument for this framing based in the psychological literature, considering the probability of one response occurring before another, where the probability of a response occurring in any given small time interval is determined by a response-specific parameter. While the argument is presented in terms of discrete time, it notes that the continuous alternative would be to consider Poisson distributions. \citet{stern1990continuum} notes that the context may be widened to that of two gamma random variables with the same shape parameter and different scale parameters, and shows that taking a shape parameter of one returns the Bradley--Terry model, whereas allowing it to tend to infinity sees the model tend to the Thurstone--Mosteller model. The idea might also be considered in the context of the discriminal process on exponential distributions of Section \ref{sec: exponential}, since the inter-arrival time of a homogeneous Poisson process with rate parameter $\lambda$ has an exponential distribution with a mean $1/\lambda$. 

More directly it is simply an expression of the standard equivalence between a multinomial distribution, in this case Bernoulli, and independent Poisson distributions conditional on their total, sometimes referred to as the ``Poisson trick'' \citep{fienberg1976log, lee2017poisson}.

\subsection{Sudden death \citep{stirzaker1999probability, vojnovic2015contest}} \label{sec: sudden death}

Consider two teams $i$ and $j$ involved in a `sudden death' shoot-out. They play a game where in each round they succeed with independent probabilities $p_i$ and $p_j$ respectively. The winner is the team who first has more successes than the other team. Let $(i \succ j)_n$ be the event that $i$ wins the `sudden death' contest in round $n$. Then
\begin{align*}
    p_{ij} &= \sum_{n=1}^{\infty}\mathbb{P}[(i \succ j)_n] \\
    &=\sum_{n=1}^{\infty}\sum_{k=0}^{n-1}p_i(1-p_j)\binom{n-1}{k}(p_ip_j)^k \\
    & \quad \quad \quad \quad \quad \quad \quad \quad ((1-p_i)(1-p_j))^{n-k-1}  \\
    &=p_i(1-p_j)\sum_{m=0}^{\infty}\sum_{k=0}^{m}\binom{m}{k}(p_ip_j)^k \\
    & \quad \quad \quad \quad \quad \quad \quad \quad ((1-p_i)(1-p_j))^{m-k} \\
    &= p_i(1-p_j)\sum_{m=0}^{\infty}(p_ip_j + (1-p_i)(1-p_j))^m \\
    &= p_i(1-p_j)\sum_{m=0}^{\infty}(2p_ip_j-p_i-p_j+1)^m \\
    &=\frac{p_i(1-p_j)}{p_i+p_j-2p_ip_j} \\
    &=\frac{p_i(1-p_j)}{p_i(1-p_j)+p_j(1-p_i)} \\
    &=\frac{\frac{p_i}{1-p_i}}{\frac{p_i}{1-p_i}+\frac{p_j}{1-p_j}} \\
    &=\frac{\pi_i}{\pi_i+\pi_j},
\end{align*}
where $\pi_i = \frac{p_i}{1-p_i}$.

Further, suppose there is an alternative contest but now the winner is the team that is the first to have $r$ more successes than the opposition. Such a contest may be understood as an aggregation of the sudden death contests described above, such that the winner is the first team to win $r$ more sudden death contests than the opposition. Based on the result above, given that there is a winner to a sudden death contest, the probability that the winner is $i$ is $q_i = p_i / (1 - p_i)$. Let $A_i$ be the event that $i$ wins and $A^{r+k}$ be the event that a result, either $i$ or $j$ winning, occurs after the winning team has won exactly $r+k$ sudden death contests then
\begin{align*}
  p_{ij} = \mathbb{P}(A_i)
            &= \sum_{k=0}^{\infty}\mathbb{P}(A_i | A^{r+k})\mathbb{P}(A^{r+k}) \\
            &= \sum_{k=0}^{\infty} \frac{q_i^{r+k} q_j^k}{q_i^{r+k} q_j^k + q_i^k q_j^{r+k}} P(A^{r+k}) \\
            &= \frac{q_i^r}{q_i^r + q_j^r} \sum_{k=0}^{\infty}P(A^{r+k}) \\
            &= \frac{q_i^r }{q_i^r + q_j^r} \\
            &= \frac{\pi_i}{\pi_i + \pi_j},
\end{align*}
where $\pi_i = q_i^r$.

\subsection{Accumulated win ratio \citep{vojnovic2015contest}} \label{sec: accumulated win ratio}

Take a sequence of matches between two players, $i$ and $j$, where the probability that team $i$ wins is proportional to the accumulated number of wins in previous matches. Suppose that the probability that $i$ wins the first match is ${\pi_i}/({\pi_i+\pi_j})$. Then consider the probability that $i$ will win the $n$th match. The claim is that this is $\pi_i / (\pi_i + \pi_j)$. We proceed to show this by induction. Define notation $(i \succ j)_n$ as meaning $i$ beats $j$ in match $n$ then our base case is
\[
\mathbb{P}[(i \succ j)_1] =\frac{\pi_i}{\pi_i+\pi_j}.
\]
Now assume the inductive hypothesis that for some $k>1$ 
\[
\mathbb{P}[(i \succ j)_k] =\frac{\pi_i}{\pi_i+\pi_j}.
\]
Then proceeding by induction
\begin{align*}
         & \mathbb{P}[(i \succ j)_{k+1}] \\
         & \quad = \mathbb{P}[(i \succ j)_{k+1}\mid(i \succ j)_k]\mathbb{P}[(i \succ j)_k] \\
         & \quad \quad \quad \quad + \mathbb{P}[(i \succ j)_{k+1}\mid(j \succ i)_k]\mathbb{P}[(j \succ i)_k] \\
         & \quad = \frac{\pi_i+1}{\pi_i+1+\pi_j}\frac{\pi_i}{\pi_i+\pi_j} + \frac{\pi_i}{\pi_i+1+\pi_j}\frac{\pi_j}{\pi_i+\pi_j} \\
         & \quad = \frac{\pi_i(\pi_i+1+\pi_j)}{(\pi_i+1+\pi_j)(\pi_i+\pi_j)} \\
         & \quad = \frac{\pi_i}{\pi_i+\pi_j} \quad .
\end{align*}

\subsection{Continuous time state transition} \label{sec: continuous time state transition}
Consider a match where the winner is the team winning at the end of a defined period of play. We choose to model the continuous state of `winning' by a continuous time Markov chain on a binary state space $I=\{i \text{ winning}, j \text{ winning}\}$,  Let the rate at which there is a switch from the state `$i \text{ winning}$'  to the state `$j \text{ winning}$' be denoted by $\pi_j$, and the rate at which the switch from the state `$j \text{ winning}$'  to the state `$i \text{ winning}$' be denoted by $\pi_i$. Then the intensity matrix is
\[
\boldsymbol{Q}=\begin{pmatrix}
  -\pi_j & \pi_j\\ 
  \pi_i & -\pi_i ,
\end{pmatrix}
\]
and the equilibrium distribution vector of this process $\boldsymbol{p}$ is such that 
\[
\boldsymbol{p}\boldsymbol{Q}=\boldsymbol{0},
\] 
and in this case is given by the probability vector $\boldsymbol{p}=\big(\frac{\pi_i}{\pi_i+\pi_j},\frac{\pi_j}{\pi_i+\pi_j}\big)$. 

Assuming that we are likely to see a large number of state changes during the course of the match or the probability of the initial state being `$i$ winning' may be approximated by $\pi_i / (\pi_i + \pi_j)$ then the probability that $i$ beats $j$ may be approximated by
\[
p_{ij} = \frac{\pi_i}{\pi_i + \pi_j}.
\]

The authors are not aware of published work that uses this continuous-time model, which might reasonably be called the ``Bradley--Terry process'' model.

%v2: The convergence to equilibrium of this Markov Chain is dictated by the number of state changes during the course of the match, which in turn will be dictated by the size of the smaller of the components of $\boldsymbol{p}$ in comparison to the length of the match. If the size of the components of $\boldsymbol{p}$ are sufficiently large, then one would expect the result to be approximately distributed according to the invariant distribution so that the probability that $i$ beats $j$ may be approximated by
%\[
%p_{ij} = \frac{\pi_i}{\pi_i + \pi_j}.
%\]

\section{Quasi-symmetry and consistent estimators} \label{sec: quasi-symmetry}

The quasi-symmetry model was proposed by \citet{caussinus1965contribution}. A matrix $C$ is quasi-symmetric if it can be decomposed such that
\[
c_{ij} = \alpha_i \beta_j \gamma_{ij},
\]
where $\gamma_{ij}=\gamma_{ji}.$ The form of this can be simplified by taking $a_{i}= \alpha_i / \beta_i$ and $s_{ij} = \beta_i \beta_j \gamma_{ij}$, so that
\[
c_{ij} = a_i s_{ij},
\]
or in matrix form
\[
C = AS,
\]
where $A$ is a diagonal matrix and $S$ is symmetric. Informally, one might think of the symmetric matrix representing the intensity of interactions, and the diagonal matrix as the relative ratings. The `gravity model' used by geographers in the study of migration is an example of a popular quasi-symmetry model.

Asymptotically, where the number of matches between each pair of teams tends to infinity and the number of teams is held constant, by the Law of Large Numbers, under a Bradley--Terry data generating process, we would expect the results matrix to be quasi-symmetric, since \[
\mathbb{E}[c_{ij}] = p_{ij}m_{ij} = \frac{\pi_i}{\pi_i + \pi_j}m_{ij} = a_{ii} s_{ij},
\]
where $s_{ij}=m_{ij}/(\pi_i+\pi_j)=s_{ji}$ and $\pi_i=a_{ii}$. So, rating methods that accord with Bradley--Terry in the case of a quasi-symmetric results matrix are consistent estimators for the Bradley--Terry model given a Bradley--Terry data generating process, and thus motivations for those rating methods are of interest in the context of this paper. This is especially so as it provides a link to a number of other, sometimes familiar, rating methods.

\subsection{PageRank \citep{daniels1969round}} \label{sec: pagerank}

\citet{daniels1969round} appears to have been the first to document the link between the Bradley--Terry model and what might now be recognized as an undamped \PageRank\ \citep{page1999pagerank}. \PageRank\ has come to be widely known as it formed the basis for the original Google search algorithm. An intuitive explanation for the way it functions is the so-called `random surfer' model. It envisages a (web-)surfer, who is randomly assigned to a node in a directed network. The random surfer then moves randomly to one of the other nodes. With a given probability they may move to any node (teleportation) or alternatively they move to a node to which there is a weighted directed edge from the node where they are currently. The probability of moving to any particular destination node if they do not teleport is set equal to the weight of the edge between the origin node and the destination node divided by the total weight of edges from the origin node. This process continues indefinitely with the proportion of time spent at each node representing the \PageRank\ for that node. What we refer to here as `undamped \PageRank' is the algorithm with the teleportation probability set to zero.

In the notation of this paper, we may take the comparison matrix to define the relevant weighted directed network, with $c_{ij}$ the weight of the directed edge from $j$ to $i$. Define $D$ as the diagonal matrix of column sums with $d_{jj} = \sum_k c_{kj}$. The undamped \PageRank\ rating vector $\boldsymbol{\alpha}_{\text{PR}}$ is the stationary distribution of the Markov chain with column-normalized comparison matrix $C D^{-1}$ as a left stochastic transition matrix. That is
\[
\boldsymbol{\alpha}_{\text{PR}} = 
C D^{-1} \boldsymbol{\alpha}_{\text{PR}} .
\]

%\footnote{By convention, it is perhaps more common to see the stationary distribution of a Markov chain written using row vectors and a right stochastic transition matrix, rather than column vectors and a left stochastic transition matrix. That is, given a transition matrix $P$, a stationary distribution would solve $\pi P = \pi$, with $\pi$ a row vector. In this thesis, $\pi$ is taken as a column vector, aiding consistency and interpretability, but the statement in this paragraph would be equivalent to taking $P = (CD^{-1})^T$ in its right stochastic form.}

While this rating is perhaps best known from its link to \PageRank, it had been previously identified as the `total influence' metric in \citet{pinski1976citation} in the context of bibliometrics. It has been independently axiomatized in \citet{altman2005ranking} and in \citet{slutzki2006scoring}. More prosaically, such a measure might be motivated in the context of sports competition by the idea of a `glory-seeker' fan, or, as \citet[p. 68]{langville2012s} terms it, the `fair weather' fan. Consider a fan who begins by selecting a team to support at random. At each step they transfer their allegiance to one of the teams that has beaten the team they previously supported. This decision is made at random in proportion to the number of their defeats that were against each team. Each team is then rated by the proportion of time that the glory-seeker has spent supporting them.

While there is a pleasing intuition to this approach, there are situations where using \PageRank\ is questionable. We present two toy examples that demonstrate just such circumstances. First, consider a five team round-robin tournament between teams A, B, C, D and E. A beats B, C and D; B beats C, D and E; C beats D and E; D beats E; and E beats A, as represented in Table \ref{tbl: tournament}.
\begin{table}[htbp!]
\caption{Five-team round-robin tournament}
\label{tbl: tournament}
\centering

\begin{tabular}{cccccc|c}

\hline
        & A & B & C & D & E & Wins \\
\hline
    A   & 0 & 1 & 1 & 1 & 0 & 3     \\
    B   & 0 & 0 & 1 & 1 & 1 & 3     \\
    C   & 0 & 0 & 0 & 1 & 1 & 2     \\
    D   & 0 & 0 & 0 & 0 & 1 & 1     \\
    E   & 1 & 0 & 0 & 0 & 0 & 1     \\
\hline
\end{tabular}

\end{table}

Undamped \PageRank\ would rate A and E joint first, because every time the glory-seeker selects team A, they will subsequently select team E, whereas standard round-robin ranking by the number of wins would rate A as joint first and E as joint last. \citet{rubinstein1980ranking} established axiomatic grounds for why number of wins should be taken as the rating in a round-robin tournament and, beyond that, it is a strong norm in competitive sport, so in this situation \PageRank\ might be deemed inappropriate.

Second, consider three teams F, G, and H. Their strengths are such that we would expect F to beat G in 2/3 of matches, F to beat H in 4/5 of matches, and G to beat H in 2/3 of matches. Now consider two tournaments between these three teams. In the first of these tournaments each team plays each other team 15 times and the proportion of results follow expectations. These results are represented in Table \ref{tbl: tournament 2}(a). In the second tournament the teams win their match-ups in the same proportions, but H plays six times more matches against both F and G; while F and G play each other the same number of times as in the first tournament, with results represented in Table \ref{tbl: tournament 2}(b).

\begin{table}[htbp!]
\caption{Three-team tournaments}
\label{tbl: tournament 2}
\centering

\subfloat[]
{
\begin{tabular}{c|ccc}
        & F & G & H  \\
\hline
    F   & 0 & 10 & 12      \\
    G   & 5 & 0 & 10      \\
    H   & 3 & 5 & 0      \\
\end{tabular}
}
\quad\quad
\subfloat[]
{
\begin{tabular}{c|ccc}
        & F & G & H  \\
\hline
    F   & 0 & 10 & 72      \\
    G   & 5 & 0 & 60      \\
    H   & 18 & 30 & 0      \\
\end{tabular}
}
\end{table}

It seems clear that based on propensity to win, in either tournament (a) or (b), team F should be ranked higher than team G and team G should be ranked higher than team H\null. \PageRank\ meets this requirement for tournament (a), but ranks H highest based on the results of tournament (b). 

In both examples, it seems that undamped \PageRank\ focuses too much on the wins of a team, ignoring the losses. In the first example, it was E's win against A that drove its high ranking rather than being balanced by its losses to B, C and D. In the second example, the number of H's wins saw it ranked highest, ignoring its higher number of losses. Therefore one suggestion to address this would be to construct a rating, $\boldsymbol{\pi} = D^{-1}\boldsymbol{\alpha}_{\text{PR}}$, by scaling the undamped \PageRank\ rating of each competitor by dividing by their number of losses.

%Let $\boldsymbol{\pi} = D^{-1}\boldsymbol{\alpha}_{\text{PR}}$. In the context of a sports ranking this then would represent the undamped PageRank rating for each team divided by the number of matches they have lost.
\[
\boldsymbol{\pi} =
D^{-1}\boldsymbol{\alpha}_{\text{PR}}=
D^{-1}CD^{-1}\boldsymbol{\alpha}_{\text{PR}}= D^{-1}C
\boldsymbol{\pi},
\]
so that $\boldsymbol{\pi}$ is an eigenvector for $\Hat{C}=D^{-1}C$. 

A vector $\boldsymbol{\pi}$ is an eigenvector for $\Hat{C}=D^{-1}C$ with an eigenvalue of 1 if and only if
\[
\sum_j c_{ij} \pi_j = d_{ii}\pi_i \quad \text{for all } i,
\]
but if $C= AS$ is quasi-symmetric such that $A$ is a diagonal matrix and $S$ is symmetric then choosing $\pi_i = a_{ii}$ yields, for all $i$,
\begin{align*}
\sum_j c_{ij} \pi_j &=
\sum_j a_{ii}s_{ij}a_{jj} \\
& = a_{ii}\sum_j s_{ji}a_{jj} \\
& = \pi_i \sum_j c_{ji} \\
& = d_{ii}\pi_i,
\end{align*}
so that the scaled undamped \PageRank\ $\boldsymbol{\pi} = D^{-1}\boldsymbol{\alpha}_{\text{PR}}$ is the diagonal component of a quasi-symmetric matrix. Equivalently it is the Bradley--Terry rating vector in the special case of a quasi-symmetric comparison matrix $C$ and thus a consistent estimator for the Bradley--Terry rating vector given a Bradley--Terry data-generating process.

In the context of bibliometrics, this rating method was proposed as the `influence weight' measure by \citet{pinski1976citation} and as `Scroogefactor' by \citet{selby2020citation}, the name we will adopt for the rating for the remainder of this section. In the bibliometric context, $c_{ij}$ within the comparison matrix represents a citation in journal $j$ of an article in journal $i$. It was motivated by noting that journals are likely to be of different sizes and that one may be interested in determining influence independent of size. The proposal was therefore to normalize the citations received by $i$ by the citations given by $i$. More recently, the `Rank Centrality' algorithm of \citet{negahban2012iterative} proposes the same estimator applied to ratio matrices, and it is also equivalent to the `Luce Spectral Ranking' of \citet{maystre2015fast} in the $k=2$ case. A more detailed discussion of these links was provided by \citet{selby2020citation, selby2024pagerank}.

As a brief illustration, we return to our examples. In the first example, with results from Table \ref{tbl: tournament}, the results do not make up a quasi-symmetric matrix, so that the Bradley--Terry rating and Scroogefactor do not align. As can be seen in Table \ref{tbl: tournament ratings}, Bradley--Terry produces the same ranking as the convention of taking the number of wins, since the vector of the number of wins is a sufficient statistic for the Bradley--Terry rating as we showed in Section \ref{sec: sufficient statistic}. Undamped \PageRank\ and Scroogefactor both rank the teams in the descending order A, B, C, D, but undamped \PageRank\ ranks E as being first equal, whereas Scroogefactor places it third. If we take number of wins to be the correct ranking, then Scroogefactor gives a more accurate ranking in placing E closer to last equal.

\begin{table}[htbp!]
\caption{Five-team round-robin tournament rating(ranking), with rating of E standardised to 1\null. \PageRank\ is undamped.}
\label{tbl: tournament ratings}
\centering

\begin{tabular}{lccccc}
\hline
        & A & B & C & D & E  \\
\hline
    Wins   & 3(1=) & 3(1=) & 2(3) & 1(4=) & 1(4=)      \\
    Bradley--Terry   & 7.57(1=) & 7.57(1=) & 2.75(3) & 1.00(4=) & 1(4=)      \\
    \PageRank   & 1.00(1=) & 0.67(3) & 0.44(4) & 0.33(5) & 1(1=)      \\
    Scroogefactor   & 3.00(1) & 2.00(2) & 0.67(4) & 0.33(5) & 1(3)  \\
\hline
\end{tabular}
\end{table}

In the second example, there is no convention such as the number of wins to anchor our methodology. But given the ratio of wins and losses for each pair, it seems clear that the teams should be ranked in descending order F, G, H. Since both results matrices are quasi-symmetric then Bradley--Terry and Scroogefactor are the same and provide a ranking in the appropriate ordering. As can be seen in Table \ref{tbl: tournament 2 rating}, this is matched by undamped \PageRank\ in the the first of the tournaments where every team plays every other the same number of times, but undamped \PageRank\ disagrees when H has a higher number of match-ups against the other two teams.

\begin{table}[htbp!]
\caption{Three-team tournament rating(ranking) with rating of H standardized to 1\null. \PageRank\ is undamped.}
\label{tbl: tournament 2 rating}
\centering

\subfloat[]
{
\begin{tabular}{lccc}
\hline
        & F & G & H  \\
\hline
    Bradley--Terry   & 4.00(1) & 2.00(2) & 1(3)      \\
    \PageRank   & 1.45(1) & 1.36(2) & 1(3)      \\
    Scroogefactor   & 4.00(1) & 2.00(2) & 1(3)      \\
\hline
\end{tabular}
}
\quad\quad
\subfloat[]
{
\begin{tabular}{lccc}
\hline
        & F & G & H  \\
\hline
     Bradley--Terry  & 4.00(1) & 2.00(2) & 1(3)      \\
     \PageRank  & 0.98(2) & 0.86(3) & 1(1)      \\
     Scroogefactor  & 4.00(1) & 2.00(2) & 1(3)      \\
\hline
\end{tabular}
}
\end{table}

\subsection{Fair bets \citep{daniels1969round}} \label{sec: fair bets}

\citet{daniels1969round} introduces an idea referred to as `fair scores'. It was elaborated on and cast in the perhaps more intuitive language of `bets' by \citet{moon1970generalized}. Both provide interesting discussions of more general approaches. More recently, \citet{slutzki2006scoring} provides an excellent summary of the approach, providing two axiomatizations for it, a presentation of a more informal motivation due to \citet{laslier1997tournament}, the link to undamped \PageRank, and an argument for why the axiomatizations may lead us to believe that the `fair bets' method is more appropriate for sports tournaments, while the undamped \PageRank\ is more suitable for citation networks.

The first of the axiomatizations shows that the `fair bets' model is the unique ranking derived under the three simultaneous requirements of uniformity, inverse proportionality to losses, and neutrality. Uniformity here requires that if a tournament outcome is balanced in the sense that every competitor has the same number of wins and losses then the competitors must be ranked equally. Inverse proportionality to losses requires that if one begins with a balanced tournament outcome, and then a single competitor's losses are multiplied by a constant then its rating will be divided by the same constant relative to the other competitors. Neutrality requires that if one begins with a balanced tournament outcome and some new matches are added between two teams where they share the wins equally then competitors will remain equally ranked. 

The second of the axiomatizations requires two axioms, consistency between a ranking and its reduced forms, and reciprocity. Reciprocity here requires that, in a two-player tournament, the ratio of the two competitors' ratings is equal to the ratio of their wins in matches between them, assuming that there are a non-zero number of matches between them. The reduced form condition considers a reduced tournament without a team $k$, with the comparison matrix modified to, in effect, reallocate results involving $k$ so that the comparison matrix is redefined as
\[
c_{ij} = 
     \begin{cases}
       0 & i=j\\
       c_{ij} + \frac{c_{ik}c_{kj}}{\sum_t c_{tk}} & \text{otherwise}. \\ 
     \end{cases}
\] 
The axiom requires that the relative ratings of two teams in any reduced tournament are equal to their ratio in the full tournament. Consistency requirements of this type are a common feature of axiomatic approaches to ranking \citep{thomson1996consistent}.

Alternatively, in-keeping with the original presentation of \citet{daniels1969round}, suppose one retrospectively wishes to assign a betting scheme to a tournament, where the loser pays to the winner an amount on the result of each match. This is subject to two conditions. First, the amount that is paid to the winner by the loser is a value dependent solely on the strength of the loser. So if $i$ beats $j$ then $i$ will receive an amount $\alpha^{\text{FB}}_j$ from $j$. Second, the betting scheme is fair. Here `fair' is taken to mean that the wagered amounts will have led to the result that betting on any team throughout the tournament will have a net gain of zero. Then one has the condition that, for all $i$,
\[
\sum_{j} c_{ij}\alpha^{\text{FB}}_j = \sum_{j} c_{ji}\alpha^{\text{FB}}_i ,
\]
where $\boldsymbol{\alpha}^{\text{FB}}$ may be taken as a rating vector for the participants, with the intuition being that one would be prepared to wager more on a strong team. 

If $C=AS$ is quasi-symmetric then 
we have for all $i$
\[
\sum_{j} a_{ii}s_{ij}\alpha^{\text{FB}}_j =
\sum_{j} a_{jj}s_{ji}\alpha^{\text{FB}}_i,
\]
so that
\[
\sum_{j} s_{ij}(a_{ii}\alpha^{\text{FB}}_j - a_{jj}\alpha^{\text{FB}}_i) = 0.
\]
Thus, $\alpha^{\text{FB}}_i = a_{ii} = \pi_i$, and the Fair Bets rating is a consistent estimator for the Bradley--Terry rating vector given a Bradley--Terry generating process.

\subsection{Wei--Kendall} \label{sec: Wei--Kendall}

The rating method introduced in \citet{wei1952algebraic} and \citet{kendall1955further} relies on an iterative application of the comparison matrix. The motivation for such a procedure might be seen by taking the five-team tournament example from Section \ref{sec: pagerank}. One might argue that ranking D and E equally is unfair as E's single victory occurred against a top-ranked team A, whereas D gained its only victory against bottom-ranked E. An approach to address this suggested by \citet{wei1952algebraic} is to weight each victory by the rating of the defeated team. The notion of inheriting the wins of a defeated opponent to inform a rating is intuitive enough that it forms the basis for the predominant rating system of the British playground game of conkers \citep{barrow2014conkers}. Under the Wei--Kendall method we would begin with a rating vector defined by the sum of wins
\[
_1\boldsymbol{\alpha}_{\text{WK}} =  {C}\boldsymbol{e} = \{3,3,2,1,1\}^T,
\]
where $\boldsymbol{e}$ is a $n \times 1$ vector of 1s. Then we assign to each team the sum of the first iteration ratings of each team they have beaten
\[
_2\boldsymbol{\alpha}_{\text{WK}} =  {C} _1\boldsymbol{\alpha}_{\text{WK}} = {C}^2\boldsymbol{e} = \{6, 4, 2, 1, 3\}^T.
\]
This second iteration measure is sometimes used in chess for tie-breaking, where it is known as the Sonneborn--Berger score 
\citep{hooper1996oxford}. But then one might reason that the victories should instead have been weighted by this updated rating. Proceeding in this way for the next five iterations we have Wei--Kendall rating vectors
\begin{align*}
    _3\boldsymbol{\alpha}_{\text{WK}} &= \{7,6,4,3,6\}^T, \\
    _4\boldsymbol{\alpha}_{\text{WK}} &= \{13,13,9,6,7\}^T, \\
    _5\boldsymbol{\alpha}_{\text{WK}} &= \{28,22,13,7,13\}^T, \\
    _6\boldsymbol{\alpha}_{\text{WK}} &= \{42,33,20,13,28\}^T, \\
    _7\boldsymbol{\alpha}_{\text{WK}} &= \{66,61,41,28,42\}^T .
\end{align*}
Note that $E$ continues to be ranked higher than $D$ and $C$.

Generalizing, one may define a series of rating vectors 
\[
_k\boldsymbol{\alpha}_{\text{WK}} = {C}^k\boldsymbol{e} .
\]
It is then natural to consider the limit, but this is clearly not convergent. However, as \citet{moon2015topics} notes, since the matrix $C$ is irreducible then by the Perron--Frobenius theorem \citep{frobenius1912matrizen} the rating vector defined by
\[
\boldsymbol{\alpha}_{\text{WK}} = \lim_{k \rightarrow \infty} \left(\frac{C}{\rho}\right)^k\boldsymbol{e} ,
\]
where $\rho$ is the dominant eigenvalue of $C$, is convergent, and this normalized limit may be thought of as a rating vector. In the case considered above this gives 
\[
\boldsymbol{\alpha}_{\text{WK}} = \{1.63, 1.38, 0.87, 0.55, 0.95\}^T .
\]

The same motivational construct can be applied to give a consistent estimator of the Bradley--Terry rating vector in the case of a Bradley--Terry data-generating process. In both cases, the idea is that we start with an intuitive rating method, based solely on the number of wins and losses for each team. It is then noted that the initial wins should not be considered equal and instead those wins should be weighted somehow to capture their relative worth. The most recent weighting vector can be used for this purpose, giving a new updated rating vector, which can itself then be used for weighting the wins. Proceeding iteratively in this way a rating is defined in the limit. 

In the case of the Wei--Kendall method, the sum of wins is used as the initial rating. Here, though, the initial rating is based on the win-loss ratio of each team, 
\[
\boldsymbol{\pi}_1 =D^{-1}C\boldsymbol{e} = \Hat{C}\boldsymbol{e}.
\]
In the re-weighting step of the Wei--Kendall method, the wins are simply weighted by the rating of the losing team. Here, the wins are likewise weighted by the rating of the losing team but then additionally, leaning on the intuition of needing to account for losses as well as wins, the resultant vector of rating-weighted wins is scaled by the losses for each team. 
\[
\boldsymbol{\pi}_k = D^{-1}C\boldsymbol{\pi}_{k-1} = \Hat{C}\boldsymbol{\pi}_{k-1}.
\]
Proceeding in this manner, we define a rating vector
\[
\boldsymbol{\pi}=\lim_{k \rightarrow \infty} \Hat{C}^k\boldsymbol{e} .
\]

Because the scaled matrix $\Hat{C}$ has unit dominant eigenvalue, by the Perron--Frobenius Theorem the limit exists and $\boldsymbol{\pi}$ is equal to the leading eigenvector of $\Hat{C}$. If additionally $\Hat{C}$ is quasi-symmetric, which it will be if $C$ is quasi-symmetric, then this leading eigenvector will be the vector of Bradley--Terry ratings. Thus, by applying the same reasoning that was used to motivate the Wei--Kendall method, we derive a consistent estimator for the Bradley--Terry rating vector given a Bradley--Terry data-generating process.

\subsection{Ratings Percentage Index}\label{sec: RPI}
A rating measure that until recently was prevalent in college sports in North America is the Ratings Percentage Index (RPI). It is commonly defined as \begin{equation*}
\begin{split}
\text{RPI} &= 25\% \times \text{Win Percent} \\
& \quad + 50\% \times \text{Opposition's Win Percent} \\
& \quad + 25\% \times \text{Opposition's Opposition's Win Percent}. \\
\end{split}
\end{equation*}
In the notation of this article, recalling that $M$ is the matrix of the number of matches, let the matrix $\Hat{M}=[\hat{m}_{ij}]$ with $\hat{m}_{ij} = m_{ij} / \sum_j m_{ij}$, so that $\hat{m}_{ij}$ is the proportion of $i$'s matches that are against team $j$. Define the win percent vector $\boldsymbol{x}=(x_1,x_2, \dots ,x_n)^T$ where $x_i=w_i/m_i = \sum_j c_{ij} / \sum_j m_{ij}$, then the RPI rating vector $\textbf{RPI}=(\text{RPI}_1,\text{RPI}_2, \dots ,\text{RPI}_n)^T$ may be defined as
\[
\textbf{RPI}=0.25  \boldsymbol{x} + 0.5  \Hat{M}\boldsymbol{x} + 0.25  \Hat{M}^2\boldsymbol{x}.
\]

An argument very much like the one in the previous section may be followed to motivate this, that we must consider the strength of opposition in aggregating wins and that we can do this iteratively. In the RPI it is assumed that the previous iterations carry information that should be retained in the overall rating and that three such applications is sufficient.

The choice of win percent as the initial rating vector and of the proportion of matches as the relevant weighting factor when taking account of the strength of opposition is intuitive, but not exclusively so. For example, one might instead take each team's win-loss ratio as the initial rating.  To account for the strength of opposition one might weight wins by the opposition's rating and then normalize those weighted wins by the number of losses, rather than weighting matches as in the RPI. The 0.25/0.5/0.25 weighting is arbitrary and indeed has been criticized as over-weighting the strength of a team's opposition and for producing perverse incentives \citep{baker2014death}. In the absence of any clear reason to do otherwise, an equal weighting might instead be applied. This would give an initial rating vector
\[
\boldsymbol{\alpha}_1 =  \Hat{C}\boldsymbol{e},
\]
and considering down to an opposition's opposition's strength as in RPI
\[
\boldsymbol{\alpha}_3 = \frac{1}{3}\Hat{C}^2\boldsymbol{\alpha}_1 + \frac{1}{3}\Hat{C}\boldsymbol{\alpha}_1 + \frac{1}{3}\boldsymbol{\alpha}_1
= \frac{1}{3}(\Hat{C}^3 + \Hat{C}^2 + \Hat{C})\boldsymbol{e} .
\]
Clearly there is no particular reason to stop after recursively considering two levels of opposition antecedents and so one might more generally consider
\[
\boldsymbol{\pi} = \lim_{r \rightarrow \infty}\frac{1}{r}\sum^{r}_{k=1}\Hat{C}^k \boldsymbol{e} .
\]
This is the row sum vector of the Cesaro average for $\Hat{C}$ and so
\[
\boldsymbol{\pi} = \lim_{k \rightarrow \infty}\Hat{C}^k \boldsymbol{e} ,
\]
giving the same result as in the previous section. And so we have that an RPI-style rating based on win-loss ratios is a consistent estimator for the Bradley--Terry rating vector, given a Bradley--Terry data-generating process.

%In fact, equal weighting is not the unique solution giving this result. \citet[ch.15]{moon2015topics} describes a model proposed by \citet{katz1953new} and \citet{thompson1958lectures}. Let a rating vector for a tournament with comparison matrix $X$ be proportional to 
%\[
%(X + aX^2 + a^2X^3 + a^3X^4 +....)\boldsymbol{e}=X(I-aX)^{-1}\boldsymbol{e},
%\]
%where $a$ is some positive constant for which the series converges, that is $a<\rho^{-1}$ where $\rho$ is the dominant eigenvalue of $X$, then \citet{thompson1958lectures} shows that the normalised relative strengths given by $X(I-aX)^{-1}\boldsymbol{e}$ are the same as those given by the Wei--Kendall method. In the present case one can take $X=\Hat{C}$ and conclude that
%\[\boldsymbol{\pi} = \lim_{r \rightarrow \infty}\sum^{r}_{k=1}a^{k-1}\Hat{C}^k \boldsymbol{e} ,
%\]
%where $a < 1$, returns the Bradley--Terry model when $C$ is quasi-symmetric, and thus may be considered a consistent estimator for the Bradley--Terry model under a Bradley--Terry data-generating process.

\subsection{``Winner stays on'' --- Barker's algorithm}\label{sec: barker}

It is a convention in some settings, for example pub pool tables, to play on the basis of ``winner stays on'', where the winner of any match continues to play the next competitor. While rarely part of an official ranking system, it is intuitive that players who spend more games as ``reigning champion'' might be considered stronger.

Suppose that one would like to design a ``winner stays on'' tournament with the property that the long-term proportion of time spent as the ``reigning champion'' is directly proportional to their strength. For a countable collection of players, let player $i$ have a specified strength of $\pi_i$. Denoting the indicator that player $i$ is the reigning champion after the $k^{th}$ game by $T_i^k$, then the design requirement can be specified as
\[
\lim_{K\rightarrow \infty} \frac{1}{K}\sum_{k=1}^K T_i^k =\frac{\pi_i}{\sum_{i=1}^n \pi_i}.
\]

To make progress with this, one must specify the selection probability for the next opponent. Suppose that the current reigning champion is player $i$, then the probability their next opponent is chosen to be player $j$ is denoted $\phi_{ij}$. Assuming all games are conditionally independent given the players involved then this construction is a Markov chain on the player identities with transition probability of switching ``reigning champion'' from player $i$ to player $j$ given by $\phi_{ij}p_{ij}$.

This setup is akin to the Markov chain Monte Carlo problem of generating samples from a probability distribution only known up to a scaling constant. Satisfying the above requirement for the tournament is equivalent to ensuring that the constructed Markov chain has an invariant distribution that is given by the (normalized) strengths.

There are many ways that the $p_{ij}$ can be specified to achieve this goal but a natural way is to invoke reversibility by designing the chain so that it satisfies the detailed balance equations. Again, there are many choices of $p_{ij}$ here but if one wishes the game outcomes to be determined directly by a ratio involving the strengths of the teams, then the natural choice would be to use the acceptance ratio from Barker’s algorithm \citet{barker1965monte},
\[
p_{ij} = \frac{\pi_i \phi_{ij}}{\pi_i \phi_{ij} + \pi_j \phi_{ji}}.
\]

This can be interpreted as a game being decided by a Bradley--Terry type probability where the player’s strength is biased for that particular game by a multiplicative factor accounting for the imbalance of symmetry for proposing that particular opponent as their next opponent. Hence, the biased strength of player $i$ is given by $\pi_i \phi_{ij}$ which is the original strength multiplied by the proposed opponent probability of choosing $j$ which is independent of the strength of player $j$.

Suppose further that the opponent proposal distribution is symmetric, i.e. $\phi_{ij}=\phi_{ji}$ for all pairs $i$ and $j$. This would be the case if the next opponent was selected uniformly at random in a finite collection of players or if there was some local standardized symmetric proposal centered about the current player’s identity. Then, the above probability that team $i$ beats team $j$ is given by
\[
p_{ij} = \frac{\pi_i}{\pi_i +\pi_j}.
\]

%\addtocontents{toc}{\protect\setcounter{tocdepth}{-1}}

\section{Discussion}\label{sec: discussion}
Faced with ``the many routes to the ubiquitous Bradley--Terry model'', two natural questions to ask are: how are these motivations linked? and how is a recognition of these diverse motivations useful in statistical modeling?

In seeking to address the linkages, we discuss the Bradley--Terry model in the context of an exponential family of distributions \citep{darmois1935lois, pitman1936sufficient, koopman1936distributions}. This provides a direct link between perhaps the most substantial motivations, those of Sections \ref{sec: axiomatic} and Sections \ref{sec: maxent} and \ref{sec: maxlikelihood}, by showing that the motivations are the specific expression in the Bradley--Terry context of general features of exponential family models. It also provides a direct link to the motivations of Sections \ref{sec: geometric}, \ref{sec: Mallows} and \ref{sec: von Mises--Fisher}, as these are explicitly exponential family models from alternative contexts translated to be applicable to the context under consideration with the Bradley--Terry model. 

The usefulness of being able to compare motivations is illustrated in two examples, where the initial motivation for using the model comes from one motivation, but by applying the insight from another motivation we are able to substantiate and advance the method.

\subsection{The Bradley--Terry model as an exponential family of distributions}\label{sec: exp family}
Following \citet{geyer2020exponentialfamilies}, a statistical model is an exponential family of distributions if it has a log-likelihood of the form
\[
l(\theta) = \langle y, \theta \rangle - k(\theta) ,
\]
where $y$ is a vector-valued canonical statistic; $\theta$ is a vector-valued canonical parameter; $\langle . , . \rangle$ represents an inner product; and $k$ is a real-valued function, the cumulant function, which is defined such that $\nabla k(\theta) = \mathbb{E}_{\theta}(Y)$ . In seeking a maximum likelihood estimate, the derivative is taken and set equal to zero
\[
0 = \nabla l(\theta) = y - \nabla k(\theta) = y - \mathbb{E}_{\theta}(Y),
\]
by the definition of the cumulant function within an exponential family.

In the model discussed here, the likelihood is
\[
\prod_{i<j} \binom{m_{ij}}{c_{ij}} p_{ij}^{c_{ij}}(1-p_{ij})^{m_{ij}-c_{ij}},
\]
so that the log-likelihood, up to a constant term, may be taken to be
\[
\frac{1}{2}\sum_{i,j}c_{ij} \log \bigg ( \frac{p_{ij}}{1-p_{ij}} \bigg ) + m_{ij} \log (1-p_{ij}),
\]
and may be rewritten in the form
\[
l(\theta)=\frac{1}{2}\sum_{i,j}c_{ij} \theta_{ij} - m_{ij} \log (1 + e^{\theta_{ij}}),
\]
where $\theta$ is the canonical parameter, a vector of length $n(n-1)$ corresponding to the directed pairwise comparisons, and with $\theta_{ij}=\log (p_{ij}/(1-p_{ij})) $; the canonical statistic vector $y$ takes scaled outcomes $c_{ij}/2$ as its elements; and the cumulant function is 
\[
k(\theta)=\sum_{i,j}m_{ij} \log (1 + e^{\theta_{ij}})/2.
\]

What \citet{geyer2007aster} refer to as an \emph{affine canonical sub-model} may be parameterized through the linear transformation \[
\theta = a + X \beta , 
\]
where $a$ is an offset vector, $X$ is a design matrix, and $\beta$ is the canonical parameter for the sub-model, giving a log-likelihood of
\[
l(\beta) = \langle X^Ty,\beta \rangle - k_{\text{SUB}}(\beta) ,
\]
where $k_{\text{SUB}}(\beta) = k(a+X\beta)$, so that this defines a new exponential family with canonical statistic vector $X^Ty$, canonical parameter vector $\beta$, and cumulant function $k_{\text{SUB}}$.

In the context of the Bradley--Terry model, one may take $a=0, \beta=\lambda$, where $\lambda$ is the vector of log-strengths $\lambda_i = \log \pi_i$, and $X$ to be the design matrix with the columns representing the $n$ participants, and the rows representing the $n(n-1)$ directed pairwise comparisons. The entry in the row corresponding to a preference for $i$ over $j$ has 1 in column $i$, $-1$ in column $j$ and zero elsewhere. This gives a log-likelihood
\begin{align*}
    l(\lambda) 
    &= \frac{1}{2}\sum_{i,j} (c_{ij}-c_{ji})\lambda_i - \frac{1}{2}\sum_{i,j} m_{ij}\log(1+e^{\lambda_i-\lambda_j}) \\
    &= \frac{1}{2}\sum_{i,j} (2c_{ij}-m_{ij})\lambda_i - \frac{1}{2}\sum_{i,j} m_{ij}\log(1+e^{\lambda_i-\lambda_j}) \\
    &= \sum_{i,j} c_{ij}\lambda_i - \frac{1}{2}\sum_{i,j} m_{ij}(\lambda_i + \log(1+e^{\lambda_i-\lambda_j})).
\end{align*}
Define a vector of wins $\boldsymbol{w}$ by $w_i = \sum_j c_{ij}$, then
\begin{align*}
    l(\lambda) 
    &= \sum_i w_i \lambda_i - \frac{1}{2}\sum_{i,j}m_{ij}(\lambda_i + \log(1 + e^{\lambda_i-\lambda_j})),
\end{align*}
defining an exponential family where the number of wins is the vector-valued canonical statistic and log-strength is the vector-valued canonical parameter. It is a feature of an exponential family of distributions that `observed equals expected', or more precisely that the observed value of the canonical statistic vector equals its expected value under the likelihood-maximizing distribution, that is to say
\[
y = \mathbb{E}_{\hat{\theta}}(Y) = \nabla k(\hat{\theta}),
\]
which under this affine canonical sub-model translates to 
\begin{align*}
w_k  
&= \frac{1}{2}\sum_j m_{kj}\left(1 + \frac{e^{\lambda_k-\lambda_j}}{1+e^{\lambda_k-\lambda_j}}\right) \\
& \quad \quad - \frac{1}{2}\sum_i m_{ik} \frac{e^{\lambda_i-\lambda_k}}{1+e^{\lambda_i-\lambda_k}} \\
&= \sum_j m_{kj} \frac{e^{\lambda_k}}{e^{\lambda_k}+e^{\lambda_j}}  \quad \text{for all }k, 
\end{align*}
noting that $p_{kj}=e^{\lambda_k}/(e^{\lambda_k}+e^{\lambda_j})$ gives what was referred to as the retrodictive criterion in Sections \ref{sec: maxent} and \ref{sec: maxlikelihood}. 

The motivations based on wins as a sufficient statistic, maximum entropy and maximum likelihood  of Sections \ref{sec: sufficient statistic}, \ref{sec: maxent}, and \ref{sec: maxlikelihood} may thus be seen as an example of a general fact about exponential families. If one starts with a canonical statistic, then the corresponding affine sub-model, if it exists, will be uniquely determined and it will be the maximum entropy and maximum likelihood model subject to the `observed equals expected' constraint on the canonical statistic. As shown in Section \ref{sec: sufficient statistic}, the requirement to take wins as a sufficient statistic leads directly to the same statistical condition as the other axiomatic motivations presented in Section \ref{sec: axiomatic}. Thus, a consideration of the Bradley--Terry model as an exponential family of distributions gives a synthesis to the axiomatic and objective function motivations. 

%One might also see the model as the result of considering a method of moments estimation that sets the canonical statistic vector equal to its observed value.

\subsection{Motivation-switching}\label{sec: motivation-switching}
In this section we present two brief examples to illustrate the usefulness of being able to consider the Bradley--Terry model from a diverse set of motivations. They are characterized by the selection of the model being based on one motivation but then justification and advancement of the methods employed being based on the consideration of other motivations. 

\subsubsection{Sports ranking}\label{sec: sports ranking}
The Bradley--Terry model is frequently employed in analysis of sports competitions. Indeed, the original work by \citet{zermelo1929berechnung} was an analysis of competitive chess. Many times the choice of the Bradley--Terry model for sports ranking may be based on its familiarity in the context, or perhaps on an informal version of the definitional simplicity motivations of Sections \ref{sec: simplicity1} and \ref{sec: simplicity2}. However, a more principled motivation for its application could rest in its status as the unique statistical pairwise comparison model for which the number of wins is a sufficient statistic (Section \ref{sec: sufficient statistic}). Taking the number of wins as the defining ranking measure in balanced sports tournaments is a strong norm and was axiomatized in \citet{rubinstein1980ranking}. It is then natural to generalize this principle by maintaining wins as a sufficient statistic to unbalanced tournaments, where competitors may play differing number of matches against differing opponents of varying strength.

This perspective also provides a natural way to extend the principle to situations where it is points rather than wins that are taken as the determining data in round-robin tournaments, allowing for result outcomes other than win/loss. Taking points as a sufficient statistic provides a principled motivation to the use of the ties model of \citet{davidson1970extending} for unbalanced tournaments in sports where the number of points on offer for a draw is half that for a win, or in employing David Firth's alt-3 model \citep{Firth2022Alt3} for soccer, where the norm is 3 points for a win and 1 for a draw.

The geometric motivation of Section \ref{sec: geometric} and the model based on permutations of Section \ref{sec: Mallows} may also be applied to extend the situations covered by ranking in a way that is consistent with these well-established sports norms. For example, in athletics --- or track and field in North American parlance --- it is common for races to be of variable size and to have different entrants at each race. If $A$, of size $n$, is the total set of competitors, let $A_k$, of size $n_k$, be the set of competitors in race $k$, and let $r_{ik}$ be the finishing position for competitor $i \in T$ in race $k$. Then we can define a result vector of length $n$ for race $k$, $\mathbf{x}_k$, with value 
\[
\frac{(n_k+1)/2 - r_{ik}}{\sqrt{n_k(n^2_k-1)/12}}
\]
if $i \in A_k$ and zero otherwise. Consistent with Sections \ref{sec: geometric} and \ref{sec: Mallows}, a rating vector $\lambda$ can then be determined by minimizing the cumulative squared Euclidean distance
\[
\sum_k d(\mathbf{x}_k, \lambda),
\]
giving a rating consistent with the Bradley--Terry model in the pairwise comparison case, but extended to allow for multiple different competitors in each contest.

%A nice feature of this conception is that it allows one to very naturally extend the Bradley--Terry model to observed results of varying size and participation, such as one finds in other sports such as athletics, motor racing or golf, if one can find a satisfactory way of representing those results on the unit sphere. 

\subsubsection{Comparative Judgment}\label{sec: Comparative Judgement}
Comparative Judgment is a form of educational assessment. It creates ratings for a set of items by having judges rank subsets of the items. These comparisons are most commonly pairwise with the Bradley--Terry model being fitted to determine the ratings. \citet{andrich1978relationships} (Section  \ref{sec: Rasch}) is often cited in that literature and so it seems reasonable to speculate that the familiarity of the Rasch model in educational assessment may be a significant reason for the model choice. But given the nature of the outcome --- the rating of academic work --- there might be a legitimate desire to be able to demonstrate the fairness of any method used. While there are not the strong norms around number of wins as a rating measure in this context like in the sports example, the idea of maximizing entropy and in that sense minimizing the assumptions in the modeling may be attractive as a justification. 

The motivations discussed here might also influence some of the practices employed in Comparative Judgment. Often the comparisons are scheduled in order to be able to produce ratings of equivalent reliability with fewer judgments than would be achieved with random scheduling. These adaptive scheduling schemes work by scheduling comparisons between items that are similar in strength so that the information from each pairwise comparison is maximized \citep{pollitt2012comparative}. The Swiss scheduling scheme, where competitors with the same, or as similar as possible, number of wins are scheduled to play each other, is a well-known example. More sophisticated approaches use an online rating that accounts for the strength of the observed comparators in order to schedule the next comparisons. These ratings could be the Bradley--Terry ratings, but their computational expense may make them unsuitable. The consistent estimators to the Bradley--Terry model discussed in Section \ref{sec: quasi-symmetry} may provide computationally faster methods for the online rating used for scheduling, even if the final rating is based on fitting the Bradley--Terry model directly, based on the fairness justification.

\section{Concluding Remarks}
In concluding, we highlight four aspects that we hope the reader may take from this work. First is a general interest in the model. Special status is accorded to models and phenomena that become apparent from a diversity of seemingly unrelated perspectives. It is in this spirit, and with a certain affection for the Bradley--Terry model, that this work was initially undertaken. Undoubtedly some of the motivations presented here carry more weight than others. Being the unique solution to maximizing entropy subject to the retrodictive criterion will be a relevant motivation in more scenarios than being a readily hypothesized model for a sudden death contest on a difference of $r$ points. Nevertheless, the number and diversity of motivations is suggestive of the applicability and attractiveness of the model, and lays the basis for its use in a wide variety of contexts. 

Second is an appreciation for the importance of model motivations. Often the motivation for using a particular model is a pragmatic one based on goodness of fit, predictive ability, computational ease or simply familiarity to the practitioner. However, there can be scenarios where a more principled motivation matters. This is likely to be the case where there are issues of fairness involved. Such scenarios are not uncommon where the output of a model is a rating, as with the examples of official sports ranking and educational assessment. The `wins as a sufficient statistic' and `maximum entropy' motivations may be particularly pertinent in those scenarios. 

Third is an appreciation for how understanding different motivations can aid in modeling practice, as illustrated with the examples of Sections \ref{sec: sports ranking} and \ref{sec: Comparative Judgement}. The setting of the Bradley--Terry model in the context of an exponential family of distributions, and the directly related motivations, may be particularly useful in advancing or expanding its application. 

Finally, we hope the work may be useful in devising material for engaging wider audiences. Some of the subject matter that the Bradley--Terry model relates to --- ratings in general, especially when applied to fields like sports --- are ones that can be of great interest to student and outside audiences, and so it is to be hoped that this work can assist in that engagement.

\begin{acks}[Acknowledgments]
The authors would like to thank Dr Martine Barons, Dr Helen Ogden, Professor Ioannis Kosmidis and two anonymous reviewers for their constructive comments that helped improved this paper.
\end{acks}

\begin{funding}
This research was partly completed during a PhD, funding for which was provided by grants A.STAA.1617.IXH and A.STAA.1819.IXH from the UK Engineering and Physical Sciences Research Council. 
\end{funding}

\bibliographystyle{imsart-nameyear.bst}
\bibliography{bibliography.bib}

\end{document}